\documentclass[a4paper,11pt,onecolumn,twoside]{article}
\usepackage{mathrsfs}
\usepackage{mathrsfs}
\usepackage{amsfonts}
\usepackage{fancyhdr}
\usepackage{amsmath,amsfonts,amssymb,graphicx,amsthm}
\allowdisplaybreaks
\begin{document}

\title{\bf Global Smooth Solutions to the Nonisothermal Compressible Fluid Models of Korteweg Type with Large Initial Data}
\author{{\bf Zhengzheng Chen}\\
School of Mathematical Sciences,
Anhui University, Hefei, 230601,  China\\[2mm]
{\bf Lin He}\\
School of Mathematics and Statistics, Wuhan University, Wuhan, 430072,  China\\[2mm]
{\bf Huijiang Zhao}\thanks{Corresponding author. E-mail:
hhjjzhao@hotmail.com}\\
School of Mathematics and Statistics, Wuhan University, Wuhan,  430072,  China}

\date{}

\vskip 0.2cm

\maketitle

\vskip 0.2cm \arraycolsep1.5pt
\newtheorem{Lemma}{Lemma}[section]
\newtheorem{Theorem}{Theorem}[section]
\newtheorem{Definition}{Definition}[section]
\newtheorem{Proposition}{Proposition}[section]
\newtheorem{Remark}{Remark}[section]
\newtheorem{Corollary}{Corollary}[section]

\begin{abstract}
The global  solutions with large initial data for the  isothermal compressible fluid models of Korteweg type have  been studied by many authors in recent years. However, little is  known of global large solutions to the  nonisothermal compressible fluid models of Korteweg type up to now. This paper is devoted to this problem, and we are concerned with the global existence of smooth and non-vacuum solutions with large initial data  to the Cauchy problem of the one-dimensional nonisothermal compressible fluid models of Korteweg type.  The case when the viscosity coefficient $\mu(\rho)=\rho^\alpha$, the  capillarity coefficient  $\kappa(\rho)=\rho^\beta$, and the heat-conductivity coefficient $\tilde{\alpha}(\theta)=\theta^\lambda$ for some parameters $\alpha,\beta,\lambda\in \mathbb{R}$ is considered. Under some assumptions on $\alpha,\beta$ and $\lambda$, we prove the global existence and time-asymptotic  behavior  of large solutions around constant  states. The proofs are given by the elementary energy method combined with the technique developed by Y. Kanel' \cite{Y. Kanel} and the maximum principle.

\bigbreak
\noindent


\end{abstract}

\section{Introduction }
\setcounter{equation}{0}
The  motion of a one-dimensional compressible nonisothermal viscous  fluid with internal capillarity  can be described by the Korteweg-type model (see \cite{J. E. Dunn-J. Serrin-1985,J. D. Van der Waals,D. J. Korteweg-1901,H. Hattori-D. Li-1996-2,B. Haspot-2009}):
\begin{eqnarray}\label{1.1}
\left\{\begin{array}{ll}
          \rho_\tau+(\rho u)_y=0,\\[2mm]
          (\rho u)_\tau+(\rho u^2+p(\rho,\theta))_y=(\mu u_y)_y+K_y,\\[2mm]
         \displaystyle \left(\rho\left(e+\frac{|u|^2}{2}\right)\right)_t+\displaystyle\left(\rho u\left(e+\frac{|u|^2}{2}\right)+p(\rho,\theta)u\right)_y=(\tilde{\alpha}\theta_y+\mu u_y+K+W)_y,
 \end{array}\right. \tau>0, y\in\mathbb{R}.
\end{eqnarray}
Here $\tau$ and $y$ represent the time variable and the spatial variable, respectively. The Korteweg tensor $K$ and the interstitial work flux $W$  are  given by
\[\left\{\begin{array}{ll}K=\kappa\rho\rho_{yy}+\rho\kappa_y\rho_y-
          \displaystyle\frac{1}{2}\kappa_\rho\rho\rho_y-\frac{\kappa}{2}\rho_y^2,\\[2mm]
          W=-\kappa\rho\rho_yu_y.
           \end{array}\right.\]
The unknown functions are the density $\rho>0$, velocity $u$, absolute temperature $\theta>0$,  pressure $p$  and the internal energy $e$ of the fluid.  $\mu,\kappa>0$ and $\tilde{\alpha}>0$ denote the viscosity coefficient,  the capillarity  coefficient and the  heat-conductivity  coefficient, respectively.  The functions $p, e$ and the coefficients $\mu,\kappa, \tilde{\alpha}$ are functions of $\rho$ and $\theta$.  Note  that when $\kappa=0$, system (\ref{1.1}) is reduced to the compressible Navier-Stokes system.

Korteweg-type models are based on an extended version of nonequilibrium
thermodynamics, which assumes that the energy of the fluid not only depends on standard variables but also on the gradient of the density.  Thus throughout this  paper, we suppose  that the  pressure $p$ and the internal energy $e$ are given by
\begin{equation}\label{1.2}
p=R\rho\theta=Av^{-\gamma}\exp\left(\frac{\gamma-1}{R}s\right),\quad e=C_v\theta+\frac{1}{2\rho}\left(\kappa(\rho,\theta)-\theta\kappa_\theta(\rho,\theta)\right)\rho_x^2,
\end{equation}
where $s$ is the entropy, $v=\frac{1}{\rho}>0$ is the specific volume, $C_v=\frac{R}{\gamma-1}$, and $\gamma>1, A, R$ are positive constants.

The global existence and large-time behavior of solutions to the compressible fluid models of Korteweg type have been studied by many authors. For small initial data, we refer to \cite{H. Hattori-D. Li-1996, H. Hattori-D. Li-1996-2} for the  global existence of  smooth solutions around constant states in Sobolev space, \cite{Y. J. Wang-Z. Tan-2011,Tan-Zhang-2014,Y.-P. Li-2012,Wang-2015,Chen-2012-1,Chen-2012-2,Chen-2012-3,Chen-2012-4} for the large-time behavior of smooth solutions in Sobolev space, and \cite{R. Danchin-B. Desjardins-2001,B. Haspot-2009} for the global existence and uniqueness of strong  solutions in  Besov space.

 For large initial data, Bresch, Desjardins, and Lin \cite{D. Bresch-B. Desjardins-C. K. Lin-2003} investigated  the global existence of weak solutions for an  isothermal fluid with a linearly density-dependent viscosity  coefficient  and a constant capillarity coefficient in a  periodic or strip domain. Later, such a result was improved by Haspot  \cite{B. Haspot-2011} to a general density-dependent viscosity coefficient and a specific type of capillarity coefficient.  Tsyganov \cite{E-Tsyganov-2008} showed the global weak solutions for an isothermal system  with the viscosity coefficient $\mu(\rho)\equiv1$  and the capillarity coefficient $\kappa(\rho)={\rho^{-5}}$ on the interval $[0,1]$.  Charve and Haspot \cite{F-Charve-B-Haspot-2011} proved  the global existence of large strong solution to the isothermal Korteweg system with $\mu(\rho)=\varepsilon\rho$ and $\kappa(\rho)=\varepsilon^2{\rho^{-1}}$ in $\mathbb{R}$. Recently, Germain and LeFloch \cite{Germain-LeFloch-2012} studied the global existence of weak solution  for the isothermal Korteweg system with general density-dependent viscosity coefficient and  capillarity coefficient in $\mathbb{R}$. Both the vacuum and non-vacuum weak solutions were obtained in \cite{Germain-LeFloch-2012}. Moreover, Chen et al. \cite{Chen-2012-5} discussed  the global existence of  smooth and non-vacumm solutions to the isothermal Korteweg system with general density-dependent viscosity and  capillarity coefficients in the whole space $\mathbb{R}$.

 Notice that the results on the large solutions of the compressible fluid models of Korteweg type \cite{D. Bresch-B. Desjardins-C. K. Lin-2003,B. Haspot-2011,E-Tsyganov-2008,F-Charve-B-Haspot-2011,Germain-LeFloch-2012,Chen-2012-4} are all concerning the isothermal system.  To the best of our knowledge, few results have been obtained for the  global large solutions to the nonisothermal compressible fluid models of Korteweg type  up to now. The main purpose of this paper is devoted to this problem, and we are concerned with the global existence of smooth, non-vacuum solutions to the Cauchy problem of system (\ref{1.1}) with large initial data  when the viscosity coefficient $\mu$ and the capillarity coefficient $\kappa$ depend only on $\rho$,  and  the  heat-conductivity coefficient $\tilde{\alpha}$ depends only on  $\theta$.  Motivated by the works \cite{Dafermos-1982,Liu-Yang-Zhao-Zou-2014,Tan-Yang-Zhao-Zou-SIMA-2013,Kawohl-1985,Jenssen-2010} on the global large solutions  for the compressible Navier-Stokes system, these density and/or temperature dependent physical coefficients will have strong influence on the solution behavior and thus lead to difficulties in analysis. In this paper, we mainly use the technique developed by Y. Kanel' \cite{Y. Kanel} and the maximum principle to overcome those difficulties.

 Now we begin to formulate our main results. Throughout this  paper, we will assume that the physical  coefficients $\mu,\kappa$ and  $\tilde{\alpha}$ are given by
\begin{equation}\label{1.3}
\mu=\mu(\rho)=\rho^\alpha, \quad \kappa=\kappa(\rho)=\rho^\beta,\quad \tilde{\alpha}=\tilde{\alpha}(\theta)=\theta^\lambda,
\end{equation}
 where $\alpha, \beta, \lambda \in\mathbb{R}$ are some constants. Then under the assumptions of (\ref{1.2}) and  (\ref{1.3}),  the system (\ref{1.1}) in the Lagrangian  coordinates can be written  as
\begin{eqnarray}\label{1.4}
\left\{\begin{array}{ll}
          v_t-u_x=0,\\[2mm]
          u_t+p(v,\theta)_x=\displaystyle\left(\frac{u_x}{v^{\alpha+1}}\right)_x+ \displaystyle\left\{-\frac{v_{xx}}{v^{\beta+5}}+\frac{\beta+5}{2}\frac{v_x^2}{v^{\beta+6}}\right\}_x,\\[2mm]
         \displaystyle C_v\theta_t+\displaystyle p(v,\theta)u_x=\displaystyle\left(\frac{\theta^\lambda\theta_x}{v}\right)_x+\frac{u_x^2}{v^{\alpha+1}},
 \end{array}\right.
\end{eqnarray}
where $t>0$ is the time variable,  $x\in\mathbb{R}$ is the Lagrangian space variable,  $v=\frac{1}{\rho}>0$ denotes the specific volume, and the pressure $p(v,\theta)=\frac{R\theta}{v}$.  System (\ref{1.4}) is imposed  with the following initial data:
\begin{equation}\label{1.5}
(v,u,\theta)(0,x)=(v_0(x), u_0(x), \theta_0(x)),\quad \lim_{x\rightarrow\pm\infty}(v_0(x),u_0(x),\theta_0(x))=(1, 0, 1).
\end{equation}

Our first theorem is concerned with the the global existence of smooth and non-vacuum  solutions to the Cauchy problem (\ref{1.4})-(\ref{1.5}), which is stated as follows.
\begin{Theorem}
Assume that the constants $\alpha, \beta$ satisfy one of the following conditions:
\[
\aligned
&(i)\quad\alpha=0,\quad \beta=-2,\\
&(ii)\quad\alpha<-2\beta-4,\quad \beta\geq-\frac{3}{2},\\
&(iii)\quad \alpha<-\beta-\frac{5}{2},\quad -2\leq\beta<-\frac{3}{2},\\
&(iv)\quad-\infty<\alpha<+\infty,\quad -3\leq\beta<-2,\\
&(v)\quad \alpha>-2\beta-5,\quad \beta<-3,
\endaligned
\]
and  the constant $\lambda$ satisfies $\lambda\geq1$. Suppose the initial data $(v_0-1, u_0, \theta_0-1)\in H^4(\mathbb{R})\times H^3(\mathbb{R})\times H^3(\mathbb{R})$, and there exist positive constants $\bar{V}, \underline{V}, \bar{\Theta}, \underline{\Theta}$ such that
\begin{equation}\label{1.6}
\underline{V}\leq v_0(x)\leq\bar{V},\quad \underline{\Theta}\leq \theta_0(x)\leq\bar{\Theta},\quad \forall\, x\in\mathbb{R}.
\end{equation}
Then the Cauchy problem (\ref{1.4})-(\ref{1.5}) admits a unique global smooth solution $(v,u,\theta)(t,x)$ satisfying for any given constant $T>0$ that
\begin{equation}\label{1.7}
C_1^{-1}(T)\leq v(t,x)\leq C_1(T),\quad C_2^{-1}(T)\leq \theta(t,x)\leq C_2(T), \quad \forall\,(t,x)\in[0,T]\times\mathbb{R},
\end{equation}
\begin{equation}\label{1.8}
\|(v-1)(t)\|^2_{H^4(\mathbb{R})}+\|(u, \theta-1)(t)\|^2_{H^3(\mathbb{R})}+\int_0^T\left(\|v_x(s)\|^2_{H^4(\mathbb{R})}+\|(u_x,\theta_x)(s)\|^2_{H^3(\mathbb{R})}\right)ds\\
\leq C_3(T),
\end{equation}
where $C_1(T), C_2(T)$ are some positive constants depending only on $T, \alpha, \beta, \lambda, \underline{V}, \bar{V}, \bar{\Theta}, \underline{\Theta}, \|v_0-1\|_{H^1(\mathbb{R})}$, $\|u_0\|_{L^2(\mathbb{R})}$ and  $\|\theta_0-1\|_{L^2(\mathbb{R})}$,  and $C_3(T)$ is a positive constant depending only on $T, \alpha, \beta, \lambda, \underline{V}, \bar{V}, \bar{\Theta}, \underline{\Theta},\\ \|v_0-1\|_{H^4(\mathbb{R})}$ and $\|(u_0, \theta_0-1)\|_{H^3(\mathbb{R})}$.
\end{Theorem}

\begin{Remark}
Some remarks on Theorem 1.1 are given as follows.
\begin{itemize}
\item[$\bullet$]The conditions (i)-(v) and $\lambda\geq1$ are used to derive the positive lower and upper bounds for the specific  volume $v(t,x)$ and the temperature $\theta(t,x)$  by using Kanel's  method \cite{Y. Kanel} and the maximum principle (see the proof of Lemma 2.2, Corollary 2.2 and Lemmas 2.6-2.11 for details).
\item[$\bullet$] Compared to the isothermal case \cite{Chen-2012-5}, a new case (i): $\alpha=0,  \beta=-2$ is obtained in Theorem 1.1 for the parameters $\alpha$ and $\beta$ by some more delicate analysis (see the proof of Lemmas 2.3, 2.6 and 2.11 for details). We believe that some similar results can also be established for the isothermal fluid by employing the method in this paper, which will be reported in a forthcoming paper by the authors.
\end{itemize}
\end{Remark}

In Theorem 1.1,  we need the parameter $\lambda$ to satisfy  $\lambda\geq1$. Our next result will remove this restriction in another setting. However, we need $\gamma-1$ to be small as a compensation. Moreover,  the large-time behavior of solutions with large initial data to the Cauchy problem (\ref{1.4})-(\ref{1.5})  will be established. Take $(v,u,s)(t,x)$ as unknown functions,  and let $\bar{s}:=\frac{R}{\gamma-1}ln\frac{R}{A}$ be the far field of the initial entropy $s_0(x)$,
i.e.,
\[\lim_{|x|\rightarrow+\infty}s_0(x)=\lim_{|x|\rightarrow+\infty}\frac{R}{\gamma-1}ln\frac{R\theta_0(x)v_0(x)^{\gamma-1}}{A}=\bar{s},
\]
then our second global existence result on the Cauchy problem (\ref{1.4})-(\ref{1.5})  can be stated as follows.
\begin{Theorem}
Let the constants $\alpha, \beta$ satisfy one of the following conditions:
\[
\aligned
&(a)\quad\alpha=0,\quad \beta=-2,\\
&(b)\quad\beta=2\alpha-3,\quad -3\leq\beta<-2,
\endaligned
\]
and  the constant $\lambda\in\mathbb{R}$. Assume that  $\|v_0-1\|_{H^4(\mathbb{R})}+\|u_0\|_{H^3(\mathbb{R})}+\|s_0-\bar{s}\|_{H^3(\mathbb{R})}$ is bounded by some constant independent of $\gamma-1$, and (\ref{1.6}) holds for some  $\gamma-1$ independent  positive constants $\bar{V}, \underline{V}, \bar{\Theta}, \underline{\Theta}$.
Then if $\gamma-1$ is sufficiently small, the Cauchy problem (\ref{1.4})-(\ref{1.5}) admits a unique global smooth solution $(v,u,\theta)(t,x)$ satisfying
\begin{equation}\label{1.9}
C_4^{-1}\leq v(t,x)\leq C_4,\quad \underline{\Theta}\leq \theta(t,x)\leq \bar{\Theta}, \quad \forall\,(t,x)\in[0,\infty)\times\mathbb{R},
\end{equation}
\begin{equation}\label{1.10}
\aligned
&\|(v-1)(t)\|^2_{H^4(\mathbb{R})}+\left\|\left(u, \frac{\theta-1}{\sqrt{\gamma-1}}\right)(t)\right\|^2_{H^3(\mathbb{R})}+\int_0^t\left(\|v_x(s)\|^2_{H^4(\mathbb{R})}+\|(u_x,\theta_x)(s)\|^2_{H^3(\mathbb{R})}\right)ds
\\
&\leq C_5\bigg(\|v_0-1\|^2_{H^4(\mathbb{R})}+\left\|\left(u_0, \frac{\theta_0-1}{\sqrt{\gamma-1}}\right)\right\|^2_{H^3(\mathbb{R})}\bigg), \quad\forall\,t>0,
\endaligned
\end{equation}
and the following time-asymptotic behavior:
\begin{equation}\label{1.10}
\lim_{t\rightarrow+\infty}\sup_{x\in\mathbb{R}}\left\{\left|\left(v(t,x)-1, u(t,x), \theta-1\right)\right|\right\}=0.
\end{equation}
Here $C_4$ is a positive constant depending only on $\lambda, \underline{V}, \bar{V}, \bar{\Theta}, \underline{\Theta}, \|v_0-1\|_{H^1(\mathbb{R})}$, $\|u_0\|_{L^2(\mathbb{R})}$ and  $\|(\theta_0-1)/\sqrt{\gamma-1}\|_{L^2(\mathbb{R})}$,  and $C_5$ is a positive constant depending only on $\lambda, \underline{V}, \bar{V}, \bar{\Theta}, \underline{\Theta},\|v_0-1\|_{H^4(\mathbb{R})}$ and $\left\|\left(u_0, (\theta_0-1)/\sqrt{\gamma-1}\right)\right\|_{H^3(\mathbb{R})}$.
\end{Theorem}

\begin{Remark}
Several remarks concerning  Theorems 1.1-1.2 are listed below.
\begin{itemize}
\item[$\bullet$] The conditions (a)-(b) and $\lambda\in\mathbb{R}$ are used to deduce the uniform-in-time lower and upper bounds for the specific  volume $v(t,x)$ and the temperature $\theta(t,x)$, and some uniform-in-time energy estimates of solutions to the Cauchy problem (\ref{1.4})-(\ref{1.5}) (see the proof of Theorem 1.2 for details).

\item[$\bullet$] In Theorem 1.2, when $\gamma-1$ is sufficiently  small, although $\left\|(\theta_0-1)/\sqrt{\gamma-1}\right\|_{H^1(\mathbb{R})}$ is small (see (\ref{3.28})), $\|(v_0-1,u_0,s_0-1)\|_{H^1(\mathbb{R})}$ can be large. Furthermore, since we only need $\gamma-1$ is sufficiently small such that $(\gamma-1)\left(\|v_0-1\|_2^2+\|(u_0,s_0-\bar{s})\|_1^2\right)<C$ for some positive constant $C$ independent of $\gamma-1$ (see (\ref{3.30}) for details), the norm $\|v_0-1\|_{H^4(\mathbb{R})}$, $\|u_0\|_{H^3(\mathbb{R})}$ and $\|\theta_{0xx}/\sqrt{\gamma-1})\|_{H^1(\mathbb{R})}$ can be  large. This is a Nishida-Smoller type result \cite{Nishida-Smoller-1973} with large initial data.

\item[$\bullet$] In Theorem 1.2, the parameter $\lambda$ can be any real constant. In fact,  Theorem 1.2 also holds when  the heat-conductivity coefficient $\tilde{\alpha}(\theta)$ is a positive function of $\theta$ with some smoothness condition. This can easily be verified from the proof of Theorem 1.2 with a slight  modification.

\item[$\bullet$] In Theorems 1.1-1.2, we obtain the global smooth solutions with large initial data for the nonisothermal compressible fluid models of Korteweg type when the physical coefficients are some power  functions of the density  $\rho$ or the temperature $\theta$.  Thus it is interesting to study some more complex and challenging cases such as when the viscosity coefficient, capillarity coefficient and heat-conductivity coefficient depend on both $\rho$ and $\theta$, which  will be  pursued by the authors  in the future.
\end{itemize}
\end{Remark}

Now we outline the main ideas to deduce the main  Theorems 1.1-1.2. Compared with  the isothermal
case [4, 8, 14, 15, 16, 19], an energy equation is involved in our problem (1.4)-(1.5), which brings
 more difficulties to prove the global existence of large solutions. The key point in the proof of Theorem 1.1 is to derive  the positive  lower and upper bounds for the specific volume $v(t,x)$ and the temperature $\theta(t,x)$. Here unlike the case of the nonisentropic compressible Navier-Stokes system \cite{Kazhikhov-Shelukhin-1977}, the argument in \cite{Kazhikhov-Shelukhin-1977} can't give the desired bounds on $v(t,x)$ and $\theta(t,x)$ for our problem (\ref{1.4})-(\ref{1.5}) due to the effect of the Korteweg tensor. To achieve these, we mainly use the method of Y. Kanel' \cite{Y. Kanel} and the maximum principle. Firstly,  an estimate of $\int_{\mathbb{R}}\frac{v_x^2}{v^{\beta+5}}dx$ appears in the basic energy estimates due to the Korteweg term (see Lemma 2.1). From this, one can deduce easily the lower and upper bounds on $v(t,x)$ for the case (iv) of Theorem 1.1 by employing  Kanel's method. Secondly, we can also perform an estimate on $\int_{\mathbb{R}}\frac{v_x^2}{v^{2\alpha+2}}dx$ from the movement equation $(\ref{1.4})_2$ under some assumptions on the parameters $\alpha$ and $\beta$ (see Lemmas 2.3-2.5 and Corollary 2.2), then by using Kanel's method again, some lower and upper bounds of $v(t,x)$ for other cases of Theorem 1.1 can be achieved in terms of $\|\theta^{1-\lambda}\|_{L^\infty([0,T_1]\times\mathbb{R})}$ (see Lemma 2.6). Thirdly, we apply
the maximum principle for the energy equation $(\ref{1.4})_3$ to obtain a lower bound
estimate on $\theta(t,x)$ in terms of $\|v^{\alpha-1}\|_{L^\infty([0,T_1]\times\mathbb{R})}$ in Lemma 2.7. The above estimates together with the upper bound estimate of $\theta(t,x)$ given in Lemma 2.10 then yield the desired lower and upper bounds of $v(t,x)$ and $\theta(t,x)$ provided that the parameters $\alpha,\beta$ and $\lambda$ satisfy  the assumptions  of Theorem 1.1 (see Corollary 2.3 and Lemma 2.11). Having obtained the  lower and upper  bounds on  $v(t,x)$ and $\theta(t,x)$, the higher order energy estimates of solutions  can be deduced easily and then Theorem 1.1 follows by the standard continuation argument. We mention that  the proof of Theorem 1.1 is greatly inspired by our previous works \cite{Chen-2012-5, Tan-Yang-Zhao-Zou-SIMA-2013}.

To prove Theorem 1.2, the main idea  is to make the a priori assumption (\ref{3.1}) for the temperature $\theta(t,x)$. Then under this a priori assumption and the conditions of Theorem 1.2, we can prove the uniform-in-time lower and upper bounds for the specific volume $v(t,x)$ and some uniform-in-time energy estimates for solutions to the Cauchy problem (\ref{1.4})-(\ref{1.5}) by using Kanel's method \cite{Y. Kanel} and some elaborate estimates. Based on these, we can close the a priori assumption (\ref{3.1}) by employing the Sobolev inequality and the smallness assumptions on $\gamma-1$. Then Theorem 1.2 follows  by combining the local existence result and the energy estimates.

 Before concluding this section,  we should point out  that there are extensive studies on the global smooth and non-vacuum solutions with large data for  the compressible Navier-Stokes system, cf. \cite{Kazhikhov-Shelukhin-1977,Y. Kanel,Dafermos-1982,Liu-Yang-Zhao-Zou-2014,Tan-Yang-Zhao-Zou-SIMA-2013,Kawohl-1985,Jenssen-2010,D. Hoff-1998,Jiang-Zhang-2003,Jiang-1999,Jiang-2002,Li-2014,Okada-Kawashima-1983} and the references therein. A recent progress in this aspect is due to Li and Liang \cite{Li-2014}, where the authors obtained the global existence and large-time behavior of smooth solutions to the one-dimensional compressible Navier-Stokes system with large initial data  and general adiabatic exponent $\gamma$. They deduce the uniform positive lower and upper bounds on  the temperature $\theta(t,x)$ by multiplying the energy equation by $(\theta-2)_+:=\max\{\theta-2,0\}$ and then some elaborate  energy analysis. However, the method in \cite{Li-2014} can't be applied to the nonisothermal compressible fluid models of Korteweg type because of the high nonlinearity of the Korteweg tensor.  The problem on how to get the large-time behavior of large solutions to the nonisothermal compressible fluid models of Korteweg type with general constant $\gamma>1$ will be left for the forthcoming future.

The rest of this paper is organized as follows.  Sections 2 and 3  are  devoted to the proof of Theorems 1.1 and 1.2, respectively.

{\bf Notations:} Throughout this paper, $C$  denotes some generic constant which may vary in
different estimates. If the dependence needs to be explicitly pointed
out, the notation $C(\cdot,\cdots,\cdot)$ or $C_i(\cdot,\cdots,\cdot)(i\in {\mathbb{N}})$ is used.  For function spaces, $H^l(\mathbb{R})$ is the
usual $l$-th order Sobolev space with its norm
\[\|f\|_{l}=\left(\sum_{i=0}^{l}\|\partial_x^if\|^2\right)^{\frac{1}{2}}\quad with\quad \|\cdot\|\triangleq\|\cdot\|_{L^2},\]
and $L^p(\mathbb{R})(1\leq p\leq+\infty)$ denotes  the usual  Lebesgue
space with the norm  $\|\cdot\|_{L^p}$. Finally,  $\|\cdot\|_{L_{T,x}^\infty}$ stands for the norm $\|\cdot\|_{L^\infty([0,T]\times\mathbb{R})}$.

\section{Proof of Theorem 1.1}
\setcounter{equation}{0}
This section is devoted to proving Theorem 1.1. To this end, we seek the solutions of the Cauchy problem (\ref{1.4})-(\ref{1.5}) in the following function space:
\[\aligned
&X(0, T; m_1, M_1; m_2, M_2)\\
&:=\left\{(\phi, \psi)(t,x)\left|
\begin{array}{c}
v(t,x)-1\in C(0, T; H^4(\mathbb{R}))\cap C^1(0, T; H^{2}(\mathbb{R}))\\[2mm]
(u(t,x),\theta(t,x)-1) \in C(0, T; (H^{3}(\mathbb{R}))^2)\cap C^1(0, T; (H^{1}(\mathbb{R}))^2)\\[2mm]
(v_x,u_x,\theta_x)(t,x)\in L^2(0, T; H^{4}(\mathbb{R})\times H^{3}(\mathbb{R})\times H^{3}(\mathbb{R}))\\[2mm]
m_1\leq v\leq M_1,\quad m_2\leq \theta\leq M_2\end{array}
\right.\right\},
\endaligned\]
where  $M_i\geq m_i>0 (i=1,2)$ and $T>0$ are some positive  constants.

Under the assumptions listed in Theorem 1.1, we can get the following local existence result.
\begin{Proposition} [Local existence]
 Under the assumptions of  Theorem 1.1, there exists a sufficiently small positive constant $t_1$ depending only on $\alpha, \beta, \lambda, \underline{V}, \bar{V},\underline{\Theta}, \bar{\Theta}$, $\|v_0-1\|_4$, $\|u_0\|_3$ and $\|\theta_0-1\|_3$ such that the Cauchy  problem (\ref{1.4})-(\ref{1.5}) admits a unique smooth solution $(v-1, u, \theta-1)(t,x)\in X(0,t_1; \frac{\underline{V}}{2}, 2\bar{V}; \frac{\underline{\Theta}}{2}, 2\bar{\Theta} )$ and
 \[\displaystyle\sup_{[0,t_1]}\{\|(v-1)(t)\|_4+\|(u, \theta-1)(t)\|_3\}\leq b\{\|v_0-1\|_4+\|(u_0, \theta_0-1)\|_3\},\]
  where $b>1$ is a positive constant depending only on  $\underline{V}, \bar{V}$.
 \end{Proposition}

The proof of Proposition 2.1 can be done  by using the dual argument and iteration technique,  which is similar to that of Theorem 2.1 in \cite{H. Hattori-D. Li-1996-2}  and  thus omitted here for brevity. To prove Theorem 1.1, one needs to show the following a priori estimates by the standard continuation argument.
\begin{Proposition} [A priori estimates]
 Under the assumptions of Theorem 1.1, suppose that $(v-1,u, \theta-1)(t,x)\in X(0, T_1; V_0, V_1; \Theta_0, \Theta_1)$ is a solution of the Cauchy problem (\ref{1.4})-(\ref{1.5}) for some  positive constants $0<T_1\leq T$ and $V_0, V_1, \Theta_0, \Theta_1 >0$.  Then there exist  positive constants $C_1(T), C_2(T)$ and $C_3(T)$ which are independent of $V_0, V_1,\Theta_0, \Theta_1 >0$  such that the following estimates hold:
\begin{equation}\label{2.1}
C_1^{-1}(T)\leq v(t,x)\leq C_1(T),\quad C_2^{-1}(T)\leq \theta(t,x)\leq C_2(T),\quad\forall\,(t,x)\in[0,T_1]\times\mathbb{R},
\end{equation}
 \begin{equation}\label{2.2}
\|(v-1)(t)\|^2_{4}+\|(u,\theta-1)(t)\|^2_{3}+\int_0^t\left(\|v_x(s)\|^2_{4}+\|(u_x, \theta_x)(s)\|^2_{3}\right)ds\leq C_3(T),\quad\forall\, t\in[0,T_1].
\end{equation}
\end{Proposition}

 In the rest of this section, we focus  on deducing Proposition 2.2.  We start with the following key lemma.
\begin{Lemma}[Basic energy estimates]
 Under the assumptions of Proposition 2.2, we have
\begin{equation}\label{2.3}
\aligned
 &\displaystyle\int_{\mathbb{R}}\eta(v,u,\theta)dx+\int_{\mathbb{R}}\frac{v_x^2}{2v^{\beta+5}}dx
 +\displaystyle\int_0^t\int_{\mathbb{R}}\left(\frac{\theta_x^2}{\theta^{2-\lambda}v}+\frac{u_x^2}{\theta v^{\alpha+1}}\right) dxd\tau\\
 &=\displaystyle\int_{\mathbb{R}}\eta(v_0,u_0,\theta_0)dx+\int_{\mathbb{R}}\frac{v_{0x}^2}{2v_0^{\beta+5}}dx
 \endaligned
 \end{equation}
for all $t\in[0, T_1]$, where the function $\eta(v,u,\theta)$ is defined by $(\ref{2.4})$.
 \end{Lemma}
\noindent{\bf Proof.}~~Let
\begin{equation}\label{2.4}
  \eta(v,u,\theta):=R\phi(v)+\frac{u^2}{2}+C_v\phi(\theta) \quad \mbox{with}\quad \phi(x)=x-lnx-1,
\end{equation}
then it is easy to see  that  $\eta(v,u,\theta)$ is a convex entropy to system (\ref{1.4}),  which  satisfies
 \begin{equation}\label{2.5}
\aligned
&\displaystyle\eta(v,u,\theta)_t+\left\{(P-R)u-\frac{uu_x}{v^{\alpha+1}}-\frac{\theta-1}{\theta}\frac{\theta^\lambda\theta_x}{v}\right\}_x
 +\displaystyle\frac{\theta^{\lambda-2}\theta_x^2}{v}+\displaystyle\frac{u_x^2}{\theta v^{\alpha+1}}\\
&=u\left\{-\frac{v_{xx}}{v^{\beta+5}}+\frac{\beta+5}{2}\frac{v_x^2}{v^{\beta+6}}\right\}_x.
\endaligned
\end{equation}
Using equation $(\ref{1.4})_1$,  we have by a direct computation that
\begin{equation}\label{2.6}
\aligned
u\left\{-\frac{v_{xx}}{v^{\beta+5}}+\frac{\beta+5}{2}\frac{v_x^2}{v^{\beta+6}}\right\}_x&=\displaystyle\{\cdots\}_x-u_x\left\{-\frac{v_{xx}}{v^{\beta+5}}+\frac{\beta+5}{2}\frac{v_x^2}{v^{\beta+6}}\right\}\\
&=\displaystyle\{\cdots\}_x-v_t\left\{-\frac{v_{xx}}{v^{\beta+5}}+\frac{\beta+5}{2}\frac{v_x^2}{v^{\beta+6}}\right\}\\
&=\displaystyle\{\cdots\}_x-\displaystyle\left(\frac{v_x^2}{2v^{\beta+5}}\right)_t.
\endaligned
\end{equation}
Here and hereafter,  $\{\cdots\}_x$ denotes the terms which will disappear after integrating with respect to $x$.

Combining (\ref{2.5}) and (\ref{2.6}), and integrating the resultant equation over $[0,t]\times\mathbb{R}$, then  we can get (\ref{2.3})  immediately. The proof of Lemma 2.1 is completed.

Based on Lemma 2.1, we now  show  the lower and upper  bounds of the specific volume $v(t,x)$ for the case (iv) of Theorem 1.1 by using Kanel's method \cite{Y. Kanel}.
\begin{Lemma} [Lower and upper  bounds of $v(t,x)$ for the case (iv) of Theorem 1.1]Let the conditions in Proposition 2.2 hold  . Suppose further that  $-3\leq\beta<-2$, then there exists a positive constant $C_6$ depending only on $ \underline{V}, \bar{V}, \underline{\Theta}, \bar{\Theta}$ and $\|(v_0-1,v_{0x}, u_0, \theta_0-1)\|$ such that
\begin{equation}\label{2.7}
C^{-1}_6\leq v(t,x)\leq C_6
 \end{equation}
for all $(t,x)\in[0, T_1]\times\mathbb{R}$.
\end{Lemma}
\noindent{\bf Proof.}~~
Define
\begin{equation}\label{2.8}
\Psi(v)=\int_1^{v}\frac{\sqrt{\phi(z)}}{z^{\frac{\beta+5}{2}}}dz,
\end{equation}
where the function $\phi(\cdot)$ is defined in (\ref{2.4}). Then  we have
\begin{equation}\label{2.9}
\Psi(v)
\longrightarrow\left\{\begin{array}{rl}
&-\infty, \quad as\,\,v\rightarrow0, \,\,if \,\,\beta\geq-3,\\[3mm]
&+\infty, \quad as\,\,v\rightarrow+\infty, \,\,if \,\,\beta<-2.
\end{array}
\right.
 \end{equation}
Furthermore,  (\ref{2.3}) together with  the Cauchy inequality implies that
\begin{equation}\label{2.10}
\begin{array}{rl}
|\Psi(v(t,x))|=&\left|\displaystyle\int_{-\infty}^x\frac{\partial}{\partial y}\Psi(v(t,y)dy\right|\\[3mm]
=&\left|\displaystyle\int_{-\infty}^x\sqrt{\phi(v)}\frac{v_y}{v^{\frac{\beta+5}{2}}}dy\right|
\leq\left\|\sqrt{\phi(v)}\right\|\left\|\frac{v_x}{v^{\frac{\beta+5}{2}}}\right\|\leq C.
\end{array}
 \end{equation}
Thus if $-3\leq\beta<-2$, then  (\ref{2.7}) follows from  (\ref{2.9}) and  (\ref{2.10}) immediately.  This completes the proof of Lemma 2.2.

\begin{Remark}
From the proof of Lemma 2.2, we see that if $\beta<-2$, then $v(t,x)$ has a uniform upper bound, i.e., $v(t,x)\leq C_6$ for all $(t,x)\in[0, T_1]\times\mathbb{R}$; and  if $\beta\geq-3$, then $v(t,x)$ has a uniform lower bound, i.e., $v(t,x)\geq C_6^{-1}$ for all $(t,x)\in[0, T_1]\times\mathbb{R}$, where $C_6$ is given in Lemma 2.2.
\end{Remark}

The following lemma gives the estimate on  $\displaystyle\int_\mathbb{R}\frac{v_x^2}{v^{2\alpha+2}}dx$.
\begin{Lemma}
 Under the assumptions of Proposition 2.2, there exist a positive constant $C=C(\underline{V}, \bar{V},\underline{\Theta}, \\ \bar{\Theta})$  and a constant $C(T)>0$ depending only on $T, \alpha,\beta,\underline{V}, \bar{V}, \underline{\Theta}, \bar{\Theta}$ and  $\|(v_0-1,v_{0x}, u_0, \theta_0-1)\|$ such that for $0\leq t\leq T_1$,
\begin{equation}\label{2.11}
\begin{array}{rl}
&\displaystyle\int_\mathbb{R}\frac{v_x^2}{v^{2\alpha+2}}dx+\displaystyle\int_0^t\int_\mathbb{R}\frac{\theta v_x^2}{v^{\alpha+1}}dxd\tau+\int_0^t\int_\mathbb{R}\left[\left(\frac{v_x}{v^{\frac{\alpha+\beta+6}{2}}}\right)_x\right]^2dxd\tau
\\[4mm]
&\leq C\|(v_0-1,v_{0x}, u_0, \theta_{0}-1)\|^2(1+\|v^{-\alpha}\|_{L^\infty_{T_1,x}}\|\theta^{1-\lambda}\|_{L^\infty_{T_1,x}})
+\displaystyle\int_0^t\int_\mathbb{R}\frac{u_x^2}{v^{\alpha+1}}dxd\tau\\[2mm]
&\quad\quad+C(T)\left\{\begin{array}{ll}
\displaystyle 0,\quad if\,\,  \alpha=f(\beta) \,\,and\,\,-5\leq\beta\leq-2, \\[3mm]
\displaystyle\left\|\frac{1}{v}\right\|^{\alpha-2\beta-5}_{L^\infty_{T_1,x}},\quad if\,\,  \alpha\geq-1 \,\,and\,\, \beta<-3, \\[3mm]
\left\|v\right\|^{-\alpha+2\beta+5}_{L^\infty_{T_1,x}},\quad if\,\,  \alpha\leq-1 \,\,and\,\,\beta\geq-2,\\[3mm]
\left\|v\right\|^{2\beta+6}_{L^\infty_{T_1,x}},\quad if\,\,  \alpha>-1 \,\,and\,\, \beta\geq-2.
\end{array}\right.
\end{array}
\end{equation}
 \end{Lemma}
\noindent{\bf Proof.}~~Using  equation $(\ref{1.4})_1$, we can rewrite  $(\ref{1.4})_2$ as
\begin{equation}\label{2.12}
\left(\displaystyle\frac{v_x}{v^{\alpha+1}}\right)_t=u_t+\displaystyle p(v,\theta)_x+\left\{\frac{v_{xx}}{v^{\beta+5}}-\frac{\beta+5}{2}\frac{v_x^2}{v^{\beta+6}}\right\}_x.
\end{equation}
Multiplying (\ref{2.12}) by $\displaystyle\frac{v_x}{v^{\alpha+1}}$ yields
\begin{equation}\label{2.13}
\displaystyle\left(\frac{v_x^2}{2v^{2\alpha+2}}-u\frac{v_x}{v^{\alpha+1}}\right)_t=\displaystyle\frac{R\theta_xv_x}{v^{\alpha+2}}-
\frac{R\theta v_x^2}{v^{\alpha+3}}+
\frac{u_x^2}{v^{\alpha+1}}+\left\{\frac{v_{xx}}{v^{\beta+5}}-\frac{\beta+5}{2}\frac{v_x^2}{v^{\beta+6}}\right\}_x\frac{v_x}{v^{\alpha+1}}
+\{\cdots\}_x.\\[3mm]
\end{equation}
Due to the formula (2.20) in \cite{Chen-2012-5}, we have
\begin{equation}\label{2.14}
\displaystyle\left\{\frac{v_{xx}}{v^{\beta+5}}-\frac{\beta+5}{2}\frac{v_x^2}{v^{\beta+6}}\right\}_x\frac{v_x}{v^{\alpha+1}}
=-\displaystyle\left[\left(\frac{v_x}{v^{\frac{\alpha+\beta+6}{2}}}\right)_x\right]^2+g(\alpha,\beta)\frac{v_x^4}{v^{\alpha+\beta+8}}+\{\cdots\}_x.
\end{equation}
where the function $g(\alpha,\beta)$ is defined by
\begin{equation}\label{2.15}
g(\alpha,\beta):=\frac{(\alpha+1)^2+(\beta+5)^2}{4}-\frac{(\beta+5)(\alpha+\beta+7)}{6}.
\end{equation}
Combining  (\ref{2.13}) and (\ref{2.14}), and integrating the resultant equation in $t$ and $x$ over $[0,t]\times\mathbb{R}$, we have
\begin{eqnarray}\label{2.16}
&&\displaystyle\frac{1}{4}\int_\mathbb{R}\frac{v_x^2}{v^{2\alpha+2}}dx+\displaystyle\int_0^t\int_\mathbb{R}\frac{R\theta v_x^2}{v^{\alpha+1}}dxd\tau+\int_0^t\int_\mathbb{R}\left[\left(\frac{v_x}{v^{\frac{\alpha+\beta+6}{2}}}\right)_x\right]^2dxd\tau
\nonumber\\
&&\leq C\displaystyle\left(\int_\mathbb{R}\frac{v_{0x}^2}{v_0^{2\alpha+2}}dx+\|u(t)\|^2+\|u_0\|^2\right)
+\int_0^t\int_\mathbb{R}\frac{u_x^2}{v^{\alpha+1}}dxd\tau\nonumber\\
&&\quad\displaystyle+\underbrace{\int_0^t\int_\mathbb{R}\frac{R\theta_xv_x}{v^{\alpha+2}}dxd\tau}_{I_1}
+\underbrace{g(\alpha,\beta)\int_0^t\int_\mathbb{R}\frac{v_x^4}{v^{\alpha+\beta+8}}dxd\tau}_{I_2},
\end{eqnarray}
where we have used the fact that
\[
\int_{\mathbb{R}}\frac{u v_x}{v^{\alpha+1}}\,dx\leq\frac{1}{4}\int_{\mathbb{R}}\frac{v_x^2}{v^{2\alpha+2}}\,dx+C\int_{\mathbb{R}}u^2\,dx.
\]
By the Cauchy inequality and (\ref{2.3}), we get
\begin{equation}\label{2.17}
\aligned
&I_1\leq\displaystyle\frac{1}{2}\int_0^t\int_\mathbb{R}\frac{R\theta v_x^2}{v^{\alpha+1}}dxd\tau+C\int_0^t\int_\mathbb{R}\frac{\theta_x^2}{\theta v^{\alpha+1}}dxd\tau\\
&\quad\leq\displaystyle\frac{1}{2}\int_0^t\int_\mathbb{R}\frac{R\theta v_x^2}{v^{\alpha+1}}dxd\tau+C\|v^{-\alpha}\|_{L^\infty_{T_1,x}}\|\theta^{1-\lambda}\|_{L^\infty_{T_1,x}}\int_0^t\int_\mathbb{R}\frac{\theta_x^2}{v\theta^{2-\lambda} }dxd\tau\\
&\quad\leq\displaystyle\frac{1}{2}\int_0^t\int_\mathbb{R}\frac{R\theta v_x^2}{v^{\alpha+1}}dxd\tau+C\|v^{-\alpha}\|_{L^\infty_{T_1,x}}\|\theta^{1-\lambda}\|_{L^\infty_{T_1,x}}.
\endaligned
\end{equation}
For the term $I_2$, if the  function $g(\alpha,\beta)=0$, i.e.,
\begin{equation}\label{2.18}
\alpha=\frac{\beta+2}{3}+\frac{1}{3}\sqrt{-2\beta^2-14\beta-20}:=f(\beta), \quad -5\leq\beta\leq-2,
\end{equation}
then $I_2=0$.

 While if $g(\alpha,\beta)\neq0$, the estimate of $I_2$ is divided into three cases. Since the estimates are the same as those in \cite{Chen-2012-5}, we cite the results directly.  According to the estimates (2.22)-(2.24) in \cite{Chen-2012-5}, we have
\begin{eqnarray}\label{2.19}
I_2\leq \eta\int_0^t \left\|\left(\frac{v_x}{v^{\frac{\alpha+\beta+6}{2}}}\right)_x(\tau)\right\|^2d\tau+ C_\eta T
\left\{\begin{array}{ll}
\displaystyle\left\|\frac{1}{v}\right\|^{\alpha-2\beta-5}_{L^\infty_{T_1,x}},\quad if\,\,  \alpha\geq-1 \,\,and\,\, \beta<-3, \\[3mm]
\left\|v\right\|^{-\alpha+2\beta+5}_{L^\infty_{T_1,x}},\quad if\,\,  \alpha\leq-1 \,\,and\,\,\beta\geq-2,\\[3mm]
\left\|v\right\|^{2\beta+6}_{L^\infty_{T_1,x}},\quad if\,\,  \alpha>-1 \,\,and\,\, \beta\geq-2.
\end{array}\right.
\end{eqnarray}
Here and hereafter, $\eta>0$ denotes a suitably small constant and $C_\eta$ is a positive constant depending on $\eta$. Putting (\ref{2.17}) and (\ref{2.19}) into (\ref{2.16}), then (\ref{2.11}) follows from Lemma 2.1 and the smallness of $\eta$. This completes the proof of Lemma 2.3.

Next, we estimate the term $\displaystyle\int_0^t\int_\mathbb{R}\frac{u_x^2}{v^{\alpha+1}}dxd\tau$ on the right hand side of (\ref{2.11}).
\begin{Lemma}  Under the assumptions of Proposition 2.2, we have
\begin{equation}\label{2.20}
\aligned
&\int_\mathbb{R}\left(u^2+\frac{v_x^2}{v^{\beta+5}}\right)dx+\displaystyle\int_0^t\int_\mathbb{R}\frac{u_x^2}{v^{\alpha+1}}dxd\tau\\
&\leq C(\underline{V}, \bar{V}, \underline{\Theta}, \bar{\Theta})\left(\|(v_0-1,v_{0x}, u_0, \theta_0-1)\|^2+\displaystyle\int_0^t\int_\mathbb{R}\frac{(\theta-1)^2}{v^{1-\alpha}}dxd\tau\right)
\endaligned
\end{equation}
for all $t\in[0,T_1]$.
\end{Lemma}
\noindent{\bf Proof.}~~Multiplying $(\ref{1.4})_2$ by $u$ and integrating  the resulting equation with respect to $t$ and
$x$ over $[0,t]\times\mathbb{R}$, we have from integration by parts that
\begin{equation}\label{2.21}
\aligned
&\int_\mathbb{R}\frac{u^2}{2}dx+\displaystyle\int_0^t\int_\mathbb{R}\frac{u_x^2}{v^{\alpha+1}}dxd\tau\\
&=-\displaystyle\int_0^t\int_\mathbb{R}up(v,\theta)_xdxd\tau+\displaystyle\int_0^t\int_\mathbb{R}u\left(-\frac{v_{xx}}{v^{\beta+5}}+\frac{\beta+5}{2}\frac{v_x^2}{v^{\beta+6}}\right)dxd\tau\\
&:=I_3+I_4.
\endaligned
\end{equation}
Using  integration by parts, we obtain
\begin{equation}\label{2.22}
\aligned
I_3&=-\displaystyle\int_0^t\int_\mathbb{R}u\left(\frac{R\theta}{v}\right)_xdxd\tau\\
&=\displaystyle\int_0^t\int_\mathbb{R}u_x\frac{R(\theta-1)}{v}dxd\tau+\displaystyle\int_0^t\int_\mathbb{R}Ru_x\left(\frac{1}{v}-1\right)dxd\tau
:=I_{3}^1+I_{3}^2.
\endaligned
\end{equation}
It follows from the basic energy estimate (\ref{2.3}) and the Cauchy inequality that
\begin{equation}\label{2.23}
I_{3}^1\leq\displaystyle\frac{1}{2}\int_0^t\int_\mathbb{R}\frac{u_x^2}{v^{\alpha+1}}dxd\tau+C\displaystyle\int_0^t\int_\mathbb{R}\frac{(\theta-1)^2}{v^{1-\alpha}}dxd\tau,
\end{equation}
\begin{equation}\label{2.24}
\aligned
I_{3}^2&=-R\displaystyle\int_0^t\int_\mathbb{R}\left(\phi(v)\right)_tdxd\tau\\
&=-R\left(\int_\mathbb{R}\phi(v)dx-\int_\mathbb{R}\phi(v_0)dx\right)\\
&\leq C(\underline{V}, \bar{V}, \underline{\Theta}, \bar{\Theta})\|(v_0-1,v_{0x}, u_0, \theta_0-1)\|^2.
\endaligned
\end{equation}
Finally, we deduce  from (\ref{2.6}) that
\begin{equation}\label{2.25}
I_{4}=-\displaystyle\int_0^t\left(\int_\mathbb{R}\frac{v_x^2}{2v^{\beta+5}}dx\right)_\tau d\tau\\
=-\int_\mathbb{R}\frac{v_x^2}{2v^{\beta+5}}dx+\int_\mathbb{R}\frac{v_{0x}^2}{2v_0^{\beta+5}}dx.
\end{equation}
Then (\ref{2.20})  follows from (\ref{2.21})-(\ref{2.25}) immediately.  This completes the proof of Lemma 2.4.

To control  the term $\displaystyle\int_0^t\int_\mathbb{R}\frac{(\theta-1)^2}{v^{1-\alpha}}dxd\tau$ on the right hand side of (\ref{2.20}), we establish the following lemma.
\begin{Lemma}  Under the assumptions of Proposition 2.2, it holds for $\lambda\neq0,-1$ that
\begin{equation}\label{2.26}
 \int_0^t\max_{x\in\mathbb{R}}|\theta(\tau,x)|^\lambda\,d\tau\leq C(T),
\end{equation}
\begin{equation}\label{2.27}
 \int_0^t\max_{x\in\mathbb{R}}|\theta(\tau,x)|^{\lambda+1}\,d\tau\leq C(T)(1+\|\theta\|_{L^{\infty}_{T_1,x}}),
\end{equation}
\begin{equation}\label{2.28}
 \int_0^t\max_{x\in\mathbb{R}}|\theta(\tau,x)|^{\lambda+1}\,d\tau\leq C(T)(1+\|v\|_{L^{\infty}_{T_1,x}}),
\end{equation}
\begin{equation}\label{2.29}
 \int_0^t\int_\mathbb{R}|\theta-1|^{2}(\tau,x)\,dxd\tau\leq C(T)(1+\|\theta^{1-\lambda}\|_{L^{\infty}_{T_1,x}}),
\end{equation}
where $C(T)$ is a positive constant depending only on $T$ and  $\|(v_0-1,v_{0x}, u_0, \theta_0-1)\|$.
\end{Lemma}
\noindent{\bf Proof.}~~(\ref{2.26})-(\ref{2.28}) can be proved by using the argument developed by Kazhikhov and Shelukhin \cite{Kazhikhov-Shelukhin-1977}, and (\ref{2.29}) can be deduced  by (\ref{2.3}) and (\ref{2.26}). The details are almost the same as those  of  Lemma 2.7 and Corollary 2.1 in \cite{Tan-Yang-Zhao-Zou-SIMA-2013}, and thus omitted here. This completes the proof of Lemma 2.5.

We  derive  from (\ref{2.29}) that
\begin{equation}\label{2.30}
\aligned
\displaystyle\int_0^t\int_\mathbb{R}\frac{(\theta-1)^2}{v^{1-\alpha}}dxd\tau&\leq\|v^{1-\alpha}\|_{L^\infty_{T_1,x}}
\displaystyle\int_0^t\int_\mathbb{R}(\theta-1)^2dxd\tau\\
&\leq C(T)\|v^{1-\alpha}\|_{L^\infty_{T_1,x}}\|\theta^{1-\lambda}\|_{L^{\infty}_{T_1,x}},
\endaligned
\end{equation}
where we have used the assumption that $\|\theta^{1-\lambda}\|_{L^{\infty}_{T_1,x}}>1$ without lose of generality.

Combining Lemmas 2.3-2.4 and  (\ref{2.30}), we have
\begin{Corollary}
Under the assumptions of Proposition 2.2, there exists a positive constant $C_{7}(T)$ depending on  $T, \alpha,\beta,\underline{V}, \bar{V}, \underline{\Theta}, \bar{\Theta}$ and  $\|(v_0-1,v_{0x}, u_0, \theta_0-1)\|$ such that  for $0\leq t\leq T_1$, it holds
\begin{equation}\label{2.31}
\begin{array}{rl}
&\displaystyle\int_\mathbb{R}\frac{v_x^2}{v^{2\alpha+2}}dx+\displaystyle\int_0^t\int_\mathbb{R}\frac{\theta v_x^2}{v^{\alpha+1}}dxd\tau+\int_0^t\int_\mathbb{R}\left[\left(\frac{v_x}{v^{\frac{\alpha+\beta+6}{2}}}\right)_x\right]^2dxd\tau
\\[4mm]
&\leq C_7(T)\left\{1+\left(\|v^{-\alpha}\|_{L^\infty_{T_1,x}}+\|v^{\alpha-1}\|_{L^\infty_{T_1,x}}\right)\|\theta^{1-\lambda}\|_{L^\infty_{T_1,x}}\right\}\\[2mm]
&\quad\quad+C_7(T)\left\{\begin{array}{ll}
\displaystyle 0,\quad if\,\,  \alpha=f(\beta) \,\,and\,\,-5\leq\beta\leq-2, \\[3mm]
\displaystyle\left\|\frac{1}{v}\right\|^{\alpha-2\beta-5}_{L^\infty_{T_1,x}},\quad if\,\,  \alpha\geq-1 \,\,and\,\, \beta<-3, \\[3mm]
\left\|v\right\|^{-\alpha+2\beta+5}_{L^\infty_{T_1,x}},\quad if\,\,  \alpha\leq-1 \,\,and\,\,\beta\geq-2,\\[3mm]
\left\|v\right\|^{2\beta+6}_{L^\infty_{T_1,x}},\quad if\,\,  \alpha>-1 \,\,and\,\, \beta\geq-2,
\end{array}\right.
\end{array}
\end{equation}
and
\begin{equation}\label{2.32}
\aligned
&\int_\mathbb{R}\left(u^2+\frac{v_x^2}{v^{\beta+5}}\right)dx+\displaystyle\int_0^t\int_\mathbb{R}\frac{u_x^2}{v^{\alpha+1}}dxd\tau\\
&\leq C_7(T)\left(\|(v_0-1,v_{0x}, u_0, \theta_0-1)\|^2+\|v^{\alpha-1}\|_{L^\infty_{T_1,x}}\|\theta^{1-\lambda}\|_{L^\infty_{T_1,x}}\right).
\endaligned
\end{equation}
\end{Corollary}

Now we use Kanel's method again to show the following lemma.
\begin{Lemma} [Lower and upper  bounds of $v(t,x)$ for the cases (i)-(iii) and (v) of Theorem 1.1]
Under the assumptions of Proposition 2.2, if  the constants $\alpha, \beta$ satisfy the condition (i), then it holds
\begin{equation}\label{2.33}
\|v\|_{L^\infty_{T_1,x}}\leq C_8(T)(1+\|\theta^{1-\lambda}\|_{L^\infty_{T_1,x}});
 \end{equation}
if the constants $\alpha, \beta$ satisfy the condition (ii) or (iii), then
\begin{equation}\label{2.34}
\|v\|_{L^\infty_{T_1,x}}\leq C_9(T)(1+\|\theta^{1-\lambda}\|_{L^\infty_{T_1,x}}^{\frac{1}{1-\alpha}});
 \end{equation}
and if the constants $\alpha, \beta$ satisfy the condition (v), then
\begin{equation}\label{2.35}
\left\|\frac{1}{v}\right\|_{L^\infty_{T_1,x}}\leq C_{10}(T)(1+\|\theta^{1-\lambda}\|_{L^\infty_{T_1,x}}^{\frac{1}{\alpha}}),
\end{equation}
where $C_i(T),i=8,9,10$ are some positive constants  depending only on  $T, \alpha,\beta,\underline{V}, \bar{V}, \underline{\Theta}, \bar{\Theta}$ and  $\|(v_0-1,v_{0x}, u_0, \theta_0-1)\|$.
\end{Lemma}
\noindent{\bf Proof.}~~ Let
\begin{equation}\label{2.36}
 \Phi(v)=\int_1^{v}\frac{\sqrt{\phi(z)}}{z^{\alpha+1}}dz,
\end{equation}
then there exist two  positive constants $A_1, A_2$ such that
\begin{equation}\label{2.37}
\left|\Phi(v)\right|\geq \left\{\begin{array}{ll} A_1v^{-\alpha}-A_2,\quad as \,\,v\rightarrow 0, \quad if\,\, \alpha>0,\\[3mm]
A_1v^{\frac{1}{2}-\alpha}-A_2,\quad as \,\,v\rightarrow +\infty, \quad if\,\, \alpha<\frac{1}{2}.
\end{array}\right.
\end{equation}

For the case (i): $\alpha=0$ and $\beta=-2$,  we have $v(t,x)\geq C_6^{-1}$ for all $(t,x)\in[0,T_1]\times\mathbb{R}$ by Remark 2.1. Now we show the upper bound of $v(t,x)$.  We deduce from $(\ref{2.37})_2$, $(\ref{2.31})_1$ with $\alpha=0, \beta=-2$, and Lemma 2.1 that
\begin{equation}\label{2.38}
\aligned
v^{\frac{1}{2}}&\leq C+C\left|\Phi(v)\right|\\
&\leq C+C\left|\displaystyle\int_{-\infty}^x\frac{\partial}{\partial y}\Phi(v(t,y))dy\right|\\
&\leq C+C\left\|\sqrt{\phi(v)}\right\|\left\|\frac{v_x}{v}\right\|\\
&\leq C+C(T)\|(v_0-1,v_{0x}, u_0, \theta_0-1)\|\left(1+
\displaystyle \|\theta^{1-\lambda}\|_{L^\infty_{T_1,x}}^{\frac{1}{2}}\right),
\endaligned
\end{equation}
which implies that
\begin{equation}\label{2.39}
\left\|v\right\|_{L^\infty_{T_1,x}}^{\frac{1}{2}}
\leq C_8(T)\left(1+\|\theta^{1-\lambda}\|_{L^\infty_{T_1,x}}^{\frac{1}{2}}\right).
\end{equation}
Then (\ref{2.33}) follows from (\ref{2.39}) immediately.

For the case (ii): $\alpha<-2\beta-4$ and $\beta\geq-\frac{3}{2}$, it holds  that $\alpha<-1$ and $v(t,x)\geq C_6^{-1}$ for all $(t,x)\in[0,T_1]\times\mathbb{R}$ by Remark 2.1. On the other hand,  it follows from $(\ref{2.37})_2$, $(\ref{2.31})_2$, Lemma 2.1 and the Young inequality  that
\begin{equation}\label{2.40}
\aligned
v^{\frac{1}{2}-\alpha}&\leq C+C\left|\Phi(v)\right|\leq C+C\left|\displaystyle\int_{-\infty}^x\frac{\partial}{\partial y}\Phi(v(t,y))dy\right|\\
&\leq C+C\left\|\sqrt{\phi(v)}\right\|\left\|\frac{v_x}{v^{\alpha+1}}\right\|\\
&\leq C+C(T)\|(v_0-1,v_{0x}, u_0, \theta_0-1)\|\left\{1+
\displaystyle \left\|v\right\|^{\frac{-\alpha+2\beta+5}{2}}_{L^\infty_{T_1,x}}\right.\\
&\quad\left.+(\|v\|_{L^\infty_{T_1,x}}^{-\frac{\alpha}{2}}
+C_6^{\frac{1-\alpha}{2}})
\|\theta^{1-\lambda}\|_{L^\infty_{T_1,x}}^{\frac{1}{2}}\right\}\\
&\leq \eta\left\|v\right\|^{\frac{1}{2}-\alpha}_{L^\infty_{T_1,x}}+C(\eta,T)\left(1+\|\theta^{1-\lambda}\|_{L^\infty_{T_1,x}}^{\frac{1-2\alpha}{2(1-\alpha)}}\right),
\endaligned
\end{equation}
provided that $\alpha<-2\beta-4$ and $\alpha<0$.
Then we can get (\ref{2.34}) from (\ref{2.40}) and the smallness of $\eta$.

For the case (iii): $\alpha<-\beta-\frac{5}{2}$ and $-2\leq\beta<-\frac{3}{2}$, it follows that  $-1<-\beta-\frac{5}{2}\leq-\frac{1}{2}$ and $v(t,x)\geq C_6^{-1}$ for all $(t,x)\in[0,T_1]\times\mathbb{R}$ by Remark 2.1. To show the upper bound of $v(t,x)$, we divide this case into the following two subcases:

 {\it Subcase $(iii)_1$}: $\alpha\leq-1$ and $-2\leq\beta<-\frac{3}{2}$. For this subcase, it holds that  $\alpha<-2\beta-4$. Then similar to the case (ii), we see that (\ref{2.34}) holds for this subcase.

 {\it Subcase $(iii)_2$}: $-1<\alpha<-\beta-\frac{5}{2}$ and $-2\leq\beta<-\frac{3}{2}$. For this subcase, similar to (\ref{2.40}), we   deduce that
\begin{equation}\label{2.41}
\aligned
v^{\frac{1}{2}-\alpha}
&\leq C(T)\left\{\|(v_0-1,v_{0x}, u_0, \theta_0-1)\|+1+
\displaystyle \left\|v\right\|^{\beta+3}_{L^\infty_{T_1,x}}\right.\\
&\quad\left.+(\|v\|_{L^\infty_{T_1,x}}^{-\frac{\alpha}{2}}
+C_6^{\frac{1-\alpha}{2}})
\|\theta^{1-\lambda}\|_{L^\infty_{T_1,x}}^{\frac{1}{2}}\right\}\\
&\leq \eta\left\|v\right\|^{\frac{1}{2}-\alpha}_{L^\infty_{T_1,x}}+C(\eta,T)\left(1+\|\theta^{1-\lambda}\|_{L^\infty_{T_1,x}}^{\frac{1-2\alpha}{2(1-\alpha)}}\right),
\endaligned
\end{equation}
provided that  $\alpha<-\beta-\frac{5}{2}$ and $\alpha<0$, then (\ref{2.41}) together with the smallness of $\eta$ leads to   (\ref{2.34}).

Finally, for the case (v): $\alpha>-2\beta-5$ and $\beta<-3$, it follows that $\alpha>1$ and  $v(t,x)\leq C_6$ for all $(t,x)\in[0,T_1]\times\mathbb{R}$ by Remark 2.1.  Moreover, similar to (\ref{2.40}), we obtain
\begin{equation}\label{2.42}
\aligned
\frac{1}{v^{\alpha}}
&\leq C(T)\left\{\|(v_0-1,v_{0x}, u_0, \theta_0-1)\|+1+
\displaystyle \left\|\frac{1}{v}\right\|^{\frac{\alpha-2\beta-5}{2}}_{L^\infty_{T_1,x}}\right.\\
&\quad\left.+\left(\left\|\frac{1}{v}\right\|_{L^\infty_{T_1,x}}^{\frac{\alpha}{2}}
+C_4^{\frac{1-\alpha}{2}}\right)
\|\theta^{1-\lambda}\|_{L^\infty_{T_1,x}}^{\frac{1}{2}}\right\}\\
&\leq \eta\left\|\frac{1}{v}\right\|^{\frac{1}{2}-\alpha}_{L^\infty_{T_1,x}}+C(\eta,T)\left(1+\|\theta^{1-\lambda}\|_{L^\infty_{T_1,x}}\right)
\endaligned
\end{equation}
provided that  $\alpha>-2\beta-5$ and $\alpha>0$.   Then  (\ref{2.35}) follows from   (\ref{2.42}) and the smallness of $\eta$ immediately. This completes the proof of Lemma 2.6.

The above Lemmas give the pointwise estimates on the specific volume $v(t,x)$.  Now we turn to estimate the temperature $\theta(t,x)$. First of all, we have the following lower bound estimate on the temperature $\theta(t,x)$.
\begin{Lemma}
Under the assumptions of Proposition 2.2, it holds that
\begin{equation}\label{2.43}
\theta(t,x)\geq C(T)\frac{1}{\left(1+\|v^{\alpha-1}\|_{L^\infty_{T_1,x}}\right)}, \quad \forall\,(t,x)\in[0,T_1]\times\mathbb{R}.
\end{equation}
\end{Lemma}
\noindent{\bf Proof.}~~ Multiplying $(\ref{1.4})_3$ by $-\frac{1}{\theta^2}$ yields
\begin{equation}\label{2.44}
\aligned
C_v\left(\frac{1}{\theta}\right)_t
&=\frac{1}{\theta^2}\left[p(v,\theta)u_x-\left(\frac{\theta^\lambda\theta_x}{v}\right)_x-\frac{u_x^2}{v^{\alpha+1}}\right]\\
&=-\frac{1}{\theta^2v^{\alpha+1}}\left(u_x-\frac{R\theta v^{\alpha}}{2}
\right)^2-\frac{2\theta_x^2\theta^{\lambda-3}}{v}+\frac{R^2v^{\alpha-1}}{4}+\left(\displaystyle\frac{\theta^\lambda }{v}\left(\frac{1}{\theta}\right)_x\right)_x.
\endaligned
\end{equation}
Let
\begin{equation}\label{2.45}
\Theta(t,x)=\frac{1}{\theta}-\frac{R^2t}{4C_v}\|v^{\alpha-1}\|_{L^\infty_{T_1,x}},
\end{equation}
then we have from (\ref{2.44}) that
 \begin{eqnarray}\label{2.46}
\left\{\begin{array}{ll}
          C_v \Theta_t\leq\left(\displaystyle\frac{\theta^\lambda \Theta_x}{v}\right)_x,\\[3mm]
          \Theta(0,x)=\frac{1}{\theta_0}\leq\frac{1}{\underline{\Theta}}.
 \end{array}\right.
\end{eqnarray}
Using the maximum principle, we obtain
\begin{equation}\label{2.47}
\Theta(t,x)\leq\frac{1}{\underline{\Theta}}, \quad \forall\,(t,x)\in[0,T_1]\times\mathbb{R},
\end{equation}
which implies (\ref{2.43}). The proof of Lemma 2.7 is completed.

For the case (iv) of Theorem 1.1, we see from Lemmas 2.2 and 2.7 that the temperature $\theta(t,x)$  has a lower  bound, i.e., there exits a positive constant $C_{11}(T)$ depending only on  $T, \underline{V}, \bar{V}, \underline{\Theta}, \bar{\Theta}$ and  $\|(v_0-1,v_{0x}, u_0, \theta_0-1)\|$ such that
\begin{equation}\label{2.48}
\theta(t,x)\geq C^{-1}_{11}(T), \quad \forall\,(t,x)\in[0,T_1]\times\mathbb{R}.
\end{equation}
Next, we  prove  the upper bound of $\theta(t,x)$ for the case (iv) of Theorem 1.1. For this purpose,  we first  show  the following corollary.
\begin{Corollary}
Under the assumptions of Proposition 2.2, if the constants $\alpha,\beta$ satisfy the condition (iv)  of Theorem 1.1, then  there exit a constant $C_{12}>0$ depending on  $T, \underline{V}, \bar{V}, \underline{\Theta}, \bar{\Theta}$ and  $\|(v_0-1,v_{0x}, u_0, \theta_0-1)\|$, and  two constants $C_{13}, C_{14}>0$ depending on  $T, \alpha,\beta,\underline{V}, \bar{V}, \underline{\Theta}, \bar{\Theta}$ and  $\|(v_0-1,v_{0x}, u_0, \theta_0-1)\|$ such that  for $0\leq t\leq T_1$,
\begin{equation}\label{2.49}
\aligned
&\left\|\left(v-1, v_x, u, \sqrt{\phi(\theta)}\right)(t)\right\|^2+\int_0^t\int_{\mathbb{R}}\left(\theta^{\lambda-2}\theta_x^2+\frac{u_x^2}{\theta}\right)dxd\tau \\&\leq C_{12}\|(v_0-1,v_{0x}, u_0, \theta_0-1)\|^2,
\endaligned
\end{equation}
\begin{equation}\label{2.50}
\left\|\left(v_x, u\right)(t)\right\|^2+\int_0^t\|u_x(\tau)\|^2d\tau\leq C_{13}(T)\left(\|(v_0-1,v_{0x}, u_0, \theta_0-1)\|^2+1\right),
\end{equation}
\begin{equation}\label{2.51}
\left\|v_x(t)\right\|^2+\int_0^t\|v_x(\tau)\|_1^2d\tau\leq C_{14}(T)\left(\|(v_0-1,v_{0x}, u_0, \theta_0-1)\|^2+1\right).
\end{equation}
\end{Corollary}
\noindent{\bf Proof.}~~(\ref{2.49}) is a direct consequence of Lemmas 2.1-2.2.  (\ref{2.50}) can be obtained from Lemma 2.2 , (\ref{2.32}), (\ref{2.48}) and the assumption that $\lambda\geq1$.

For (\ref{2.51}), notice that
\begin{equation}\label{2.52}
\left[\left(\frac{v_x}{v^{\frac{\alpha+\beta+6}{2}}}\right)_x\right]^2=\frac{v_{xx}^2}{v^{\alpha+\beta+6}}-(\alpha+\beta+6)\frac{v_x^2v_{xx}}{v^{\alpha+\beta+7}}+
\frac{(\alpha+\beta+6)^2}{4}\frac{v_x^4}{v^{\alpha+\beta+8}},
\end{equation}
then it follows from (\ref{2.31}), (\ref{2.52}) and Lemma 2.2 that
\begin{equation}\label{2.53}
\aligned
\left\|v_x(t)\right\|^2+\int_0^t\|v_x(\tau)\|_1^2d\tau&\leq C(T)\left(\|(v_0-1,v_{0x}, u_0, \theta_0-1)\|^2+\|\theta^{1-\lambda}\|_{L^\infty_{T_1,x}}\right)\\
&\quad+C(T)\int_0^t\int_\mathbb{R}v_x^4dxd\tau.
\endaligned
\end{equation}
Using the Sobolev inequality,  the Cauchy inequality and (\ref{2.49}), we obtain
\begin{equation}\label{2.54}
\aligned
\int_0^t\int_\mathbb{R}v_x^4dxd\tau&\leq \int_0^t\|v_x(\tau)\|^3\|v_{xx}(\tau)\|d\tau\\
&\leq \eta \int_0^t\|v_{xx}(\tau)\|^2d\tau+C_\eta\int_0^t\sup_{0\leq\tau\leq t}\left\{\|v_x(\tau)\|^6\right\}d\tau\\
&\leq \eta \int_0^t\|v_{xx}(\tau)\|^2d\tau+C(\eta, T).
\endaligned
\end{equation}
Combining (\ref{2.53})-(\ref{2.54}), then we can get (\ref{2.51}) by (\ref{2.48}), the assumption that $\lambda\geq1$ and the smallness of $\eta$. This completes the proof of Corollary 2.2.

For the  estimate on $\|u_x(t)\|$, we have
\begin{Lemma}
Under the assumptions of Proposition 2.2, if the constants $\alpha,\beta$ satisfy the condition (iv)  of Theorem 1.1, then there exits a constant $C_{15}(T)>0$  depending only on  $T, \alpha,\beta,\lambda,\underline{V}, \bar{V}, \underline{\Theta}, \bar{\Theta}$, $\|v_0-1\|_2$, $\|u_0\|_1$ and $\|\theta_0-1\|$  such that  for $0\leq t\leq T_1$,
\begin{equation}\label{2.55}
\left\|(u_x,v_{xx})(t)\right\|^2+\int_0^t\|u_{xx}(\tau)\|^2d\tau\leq C_{15}(T)\left(1+\int_0^t\|\theta_x(\tau)\|^2d\tau\right).
\end{equation}
\end{Lemma}
\noindent{\bf Proof.}~~Multiplying $(\ref{1.4})_2$ by $-u_{xx}$ and using equation $(\ref{1.4})_1$, we have
\begin{equation}\label{2.56}
\begin{array}{rl}
 \displaystyle\left(\frac{u_x^2}{2}+\frac{v_{xx}^2}{2v^{\beta+5}}\right)_t+\frac{u_{xx}^2}{v^{\alpha+1}}
= &\displaystyle R\left(\frac{\theta_x}{v}-\frac{\theta v_x}{v^2}\right)u_{xx}
 \displaystyle+(\alpha+1)\frac{u_xv_xu_{xx}}{v^{\alpha+2}}-\displaystyle(\beta+5)\frac{v_{xx}^2u_x}{2v^{\beta+6}}
 \\[3mm]&+(\beta+5)\displaystyle\frac{v_xv_{xx}u_{xx}}{v^{\beta+6}}-
 \displaystyle(\beta+5)(\beta+6)\frac{v_{x}^3u_{xx}}{2v^{\beta+7}}
 +\{\cdots\}_x.
 \end{array}
\end{equation}
Integrating (\ref{2.56}) with respect to $t$ and $x$ over $[0,t]\times\mathbb{R}$, we have from Lemma  2.2 that
\begin{equation}\label{2.57}
 \displaystyle\|(u_x, v_{xx})(t)\|^2+\displaystyle\int_0^t\|u_{xx}(\tau)\|^2\,d\tau
 \leq C(T)\left(\|(v_{0xx},u_{0x})\|^2+I_5+I_6\right),
 \end{equation}
where
\[
\begin{array}{rl}
 &I_5=\displaystyle\int_0^t\int_\mathbb{R}\left\{\left|v_{xx}^2u_x\right|+|\theta_xu_{xx}|+|\theta v_xu_{xx}|\right\}dxd\tau,\\[3mm]
 &I_6=\displaystyle\int_0^t\int_\mathbb{R}\left\{\left|v_xv_{xx}u_{xx}\right|+\left|u_xv_xu_{xx}\right|+\left|v_{x}^3u_{xx}\right|\right\}dxd\tau.
 \end{array}
\]
It follows from  the Sobolev inequality,  the Cauchy inequality and Corollary 2.1 that
\begin{equation}\label{2.58}
\aligned
 I_5&\leq \int_0^t\left\{\|u_{x}(\tau)\|^{\frac{1}{2}}\|u_{xx}(\tau)\|^{\frac{1}{2}}\|v_{xx}(\tau)\|^2
 +\|u_{xx}(\tau)\|\|\theta_{x}(\tau)\|\right\}\,d\tau\\
 &\quad+\eta\int_0^t\|u_{xx}(\tau)\|^{2}d\tau+
 C_\eta\int_0^t\|\theta(\tau)\|^2_{L^\infty}\sup_{0\leq\tau\leq t}\{\|v_{x}(\tau)\|\}^2d\tau\\
 &\leq3\eta\int_0^t\|u_{xx}(\tau)\|^{2}d\tau+C_\eta\int_0^t\left\{\|u_{x}(\tau)\|^{\frac{2}{3}}\|v_{xx}(\tau)\|^{\frac{8}{3}}
 +\|\theta_x(\tau)\|^2+\|\theta(\tau)\|^2_{L^\infty}\right\}d\tau\\
 &\leq3\eta\int_0^t\|u_{xx}(\tau)\|^{2}d\tau+
 C(\eta,T)\left\{\int_0^t\left(\|v_{xx}(\tau)\|^{4}+\|\theta_x(\tau)\|^2\right)d\tau+1\right\},
\endaligned
\end{equation}
where we have used the fact that
\begin{equation}\label{2.59}
\aligned
  \int_0^t\|\theta(\tau)\|^2_{L^\infty}\,d\tau&\leq
 \int_0^t\sup_{x\in\mathbb{R}}\left(\frac{1}{\theta(\tau,x)}\right)^{\lambda-1}\sup_{x\in\mathbb{R}}\left|\theta(\tau,x)\right|^{\lambda+1}d\tau\\
 &\leq (C_{11}(T))^{\lambda-1}\int_0^t\sup_{x\in\mathbb{R}}\left|\theta(\tau,x)\right|^{\lambda+1}d\tau\\
 &\leq C(T)(1+\|v\|_{L^\infty_{T_1,x}})\leq C(T),
\endaligned
\end{equation}
due to lemmas 2.2 and 2.5, (\ref{2.48}) and the assumption that $\lambda\geq1$.

Similar to (\ref{2.58}), $I_6$ can be controlled by
\begin{eqnarray}\label{2.60}
I_6 &&\leq\eta\int_0^t\|u_{xx}(\tau)\|^{2}d\tau+C_\eta\int_0^t\int_{\mathbb{R}}\left(|u_xv_x|^2+|v_xv_{xx}|^2+|v_x|^6\right)dxd\tau\nonumber\\
 &&\leq\eta\int_0^t\|u_{xx}(\tau)\|^{2}d\tau+C_\eta\int_0^t\left(\|u_x(\tau)\|\|u_{xx}(\tau)\|\|v_x(\tau)\|^2\right.\nonumber\\
 &&\left.\qquad+\|v_{x}(\tau)\|\|v_{xx}(\tau)\|^3
 +\|v_{x}(\tau)\|^4\|v_{xx}(\tau)\|^2\right)d\tau\\
 &&\leq2\eta\int_0^t\|u_{xx}(\tau)\|^{2}d\tau+C_\eta\int_0^t\left(\|v_{x}(\tau)\|_1^2+1\right)\|v_{xx}(\tau)\|^2d\tau+C(\eta,T).\nonumber
\end{eqnarray}

Substituting (\ref{2.58}) and (\ref{2.60}) into (\ref{2.57}), then  (\ref{2.55}) follows from Gronwall's inequality,  Corollary 2.2, and the smallness of $\eta$.  This completes the proof of Lemma 2.8.

The next lemma gives the estimate on $\displaystyle\int_0^t\|\theta_x(\tau)\|^2d\tau$.
\begin{Lemma}
Under the assumptions of Proposition 2.2, if the constants $\alpha,\beta$ satisfy the condition (iv) of Theorem 1.1, then there exits a constant $C_{16}(T)>0$  depending only on  $T, \alpha,\beta,\lambda, \underline{V}, \bar{V}, \underline{\Theta}, \bar{\Theta}$, $\|v_0-1\|_2$, $\|u_0\|_1$, and $\|\theta_0-1\|$  such that  for $0\leq t\leq T_1$,
\begin{equation}\label{2.61}
\left\|(\theta-1)(t)\right\|^2+\int_0^t\|\theta_{x}(\tau)\|^2d\tau\leq C_{16}(T).
\end{equation}
\end{Lemma}
\noindent{\bf Proof.}~~Multiplying  $(\ref{1.4})_3$ by $\theta-1$ and integrating the resulting equation with respect to $t$ and $x$ over $[0,t]\times\mathbb{R}$, we have
\begin{equation}\label{2.62}
\aligned
\frac{C_v}{2}\int_{\mathbb{R}}(\theta-1)^2dx+\int_0^t\int_{\mathbb{R}}\frac{\theta^\lambda\theta_x^2}{v}dxd\tau=&\frac{C_v}{2}\int_{\mathbb{R}}(\theta_0-1)^2dx
\underbrace{-\int_0^t\int_{\mathbb{R}}pu_x(\theta-1)dxd\tau}_{I_7}\\
&+\underbrace{\int_0^t\int_{\mathbb{R}}\frac{u_x^2}{v^{\alpha+1}}(\theta-1)dxd\tau}_{I_8}.
\endaligned
\end{equation}
 It follows from Lemma 2.2, (\ref{2.48}) and the assumption $\lambda\geq1$ that
 \begin{equation}\label{2.63}
\int_0^t\int_{\mathbb{R}}\frac{\theta^\lambda\theta_x^2}{v}dxd\tau\geq\displaystyle \frac{(\inf_{x\in\mathbb{R}}\theta)^\lambda}{\sup_{x\in\mathbb{R}}v}\int_0^t\|\theta_x(\tau)\|^2d\tau\geq C(T)\int_0^t\|\theta_x(\tau)\|^2d\tau.
\end{equation}
Utilizing  the Cauchy inequality, the Sobolev inequality,  Corollary 2.2 and (\ref{2.59}), we get
\begin{equation}\label{2.64}
\aligned
 I_7&=\int_0^t\int_{\mathbb{R}}\left(\frac{R\theta(\theta-1)}{v}\right)_xudxd\tau\\
 &=R\int_0^t\int_{\mathbb{R}}\left(\frac{(2\theta-1)\theta_x}{v}-\frac{\theta(\theta-1)v_x}{v^2}\right)udxd\tau\\
 &\leq\eta\int_0^t\|\theta_{x}(\tau)\|^{2}d\tau+C_\eta\int_0^t\int_{\mathbb{R}}|2\theta-1|^2|u|^2dxd\tau
 +C_\eta\int_0^t\int_{\mathbb{R}}|\theta(\theta-1)uv_x|dxd\tau\\
 &\leq\eta\int_0^t\|\theta_{x}(\tau)\|^{2}d\tau+C_\eta\int_0^t\left(\|\theta(\tau)\|_{L^\infty}^2+1\right)\sup_{0\leq\tau\leq t}\{\|(u,v_x)(\tau)\|^2\}d\tau\\
 &\leq\eta\int_0^t\|\theta_{x}(\tau)\|^{2}d\tau+C(\eta,T),
 \endaligned
\end{equation}
\begin{eqnarray}\label{2.65}
 I_8&&=-\int_0^t\int_{\mathbb{R}}uu_x\left(\frac{\theta-1}{v^{\alpha+1}}\right)_xdxd\tau-
 \int_0^t\int_{\mathbb{R}}uu_{xx}\frac{\theta-1}{v^{\alpha+1}}dxd\tau\nonumber\\
 &&=-\int_0^t\int_{\mathbb{R}}\frac{uu_x\theta_x}{v^{\alpha+1}}dxd\tau
 +(\alpha+1)\int_0^t\int_{\mathbb{R}}\frac{uu_xv_x(\theta-1)}{v^{\alpha+2}}dxd\tau
 -\int_0^t\int_{\mathbb{R}}\frac{uu_{xx}(\theta-1)}{v^{\alpha+1}}dxd\tau\nonumber\\
 &&\leq\eta\int_0^t\|(\theta_x,u_{xx})(\tau)\|^2d\tau+C_\eta\int_0^t\int_{\mathbb{R}}\left(u^2u^2_{x}+u_x^2v_x^2+u^2(\theta-1)^2\right)dxd\tau\nonumber\\
 &&\leq\eta\int_0^t\|(\theta_x,u_{xx})(\tau)\|^2d\tau+C_\eta\int_0^t\|u_x(\tau)\|\|u_{xx}(\tau)\|\sup_{0\leq\tau\leq t}\{\|u(\tau)\|+\|v_x(\tau)\|^2\}d\tau\nonumber\\
 &&\qquad+C_\eta\int_0^t\left(\|\theta(\tau)\|^2_{L^\infty}+1\right)\sup_{0\leq\tau\leq t}\{\|u(\tau)\|^2\}d\tau\nonumber\\
 &&\leq2\eta\int_0^t\|(\theta_x,u_{xx})(\tau)\|^2d\tau+C_\eta\int_0^t\left(\|u_x(\tau)\|^2+\|\theta(\tau)\|^2_{L^\infty}\right)d\tau+C(\eta,T)\nonumber\\
 &&\leq2\eta\int_0^t\|(\theta_x,u_{xx})(\tau)\|^2d\tau+C(\eta,T).
\end{eqnarray}
Combining (\ref{2.62})-(\ref{2.65}), then (\ref{2.61}) follows from Lemma 2.8 and the smallness of $\eta$. This completes the proof of Lemma 2.9.

Now we give the  upper bound of $\theta(t,x)$.
\begin{Lemma}
Under the assumptions of Proposition 2.2, if the constants $\alpha,\beta$ satisfy the condition (iv) of Theorem 1.1, then there exits a constant $C_{17}(T)>0$  depending only on  $T, \alpha,\beta,\lambda,\underline{V}, \bar{V}, \underline{\Theta}, \bar{\Theta}$, $\|v_0-1\|_2$, $\|u_0\|_2$ and $\|\theta_0-1\|$  such that
\begin{equation}\label{2.66}
\left\|\theta\right\|_{L^{\infty}_{T_1,x}}\leq C_{17}(T).
\end{equation}
\end{Lemma}
\noindent{\bf Proof.}~For any $q>1$,  multiplying $(\ref{1.4})_3$ by $2q(\theta-1)^{2q-1}$, and then integrating  the
resultant equation with respect to $t$ and $x$ over $[0,t]\times\mathbb{R}$,  we have
\begin{equation}\label{2.67}
\aligned
&C_v\|(\theta-1)(t)\|^{2q}_{L^{2q}}+2q(2q-1)\int_0^t\int_{\mathbb{R}}(\theta-1)^{2p-1}\frac{\theta^\lambda\theta_x^2}{v}dxd\tau\\
&=C_v\|(\theta_0-1)(t)\|^{2q}_{L^{2q}}
\underbrace{-2qR\int_0^t\int_{\mathbb{R}}\frac{\theta u_x}{v}(\theta-1)^{2q-1}dxd\tau}_{I_9}+\underbrace{2q\int_0^t\int_{\mathbb{R}}\frac{u_x^2}{v^{\alpha+1}}(\theta-1)^{2q-1}dxd\tau}_{I_{10}}.
\endaligned
\end{equation}
The H\"{o}lder inequality implies that
\begin{equation}\label{2.68}
\aligned
|I_9+I_{10}|&\leq C\int_0^t\int_{\mathbb{R}}\left(|\theta u_x|+u_x^2\right)|\theta-1|^{2q-1}dxd\tau\\
&\leq C\int_0^t\|(\theta-1)(\tau)\|^{2q-1}_{L^{2q}}\left(\|\theta u_x(\tau)\|_{L^{2q}}+\|u^2_x(\tau)\|_{L^{2q}}\right)d\tau
\\
&\leq C\sup_{0\leq\tau\leq T_1}\{\|(\theta-1)(\tau)\|^{2q-1}_{L^{2q}}\}\int_0^t\left(\|\theta u_x(\tau)\|_{L^{2q}}+\|u^2_x(\tau)\|_{L^{2q}}\right)d\tau
\endaligned
\end{equation}
Combining (\ref{2.67})-(\ref{2.68}) yields
\begin{equation}\label{2.69}
\aligned
\sup_{0\leq\tau\leq T_1}\{\|(\theta-1)(\tau)\|_{L^{2q}}\}\leq C\int_0^t\left(\|\theta u_x(\tau)\|_{L^{2q}}+\|u^2_x(\tau)\|_{L^{2q}}\right)d\tau+C.
\endaligned
\end{equation}
Letting $q\rightarrow+\infty$, we deduce  from (\ref{2.69}), (\ref{2.50}), (\ref{2.59}), the  Sobolev inequality and Lemmas 2.8-2.9 that
\begin{eqnarray}\label{2.70}
\|\theta\|_{L^{\infty}_{T_1,x}}&&\leq C\int_0^t\left(\|\theta u_x(\tau)\|_{L^\infty}+\|u^2_x(\tau)\|_{L^\infty}\right)d\tau+C\nonumber\\
&&\leq C\int_0^t\left(\|u_x(\tau)\|\|u_{xx}(\tau)\|+\|\theta(\tau)\|_{L^\infty}^2\right)d\tau+C\nonumber\\
&&\leq C\int_0^t\left(\|u_x(\tau)\|_1^2+\|\theta(\tau)\|_{L^\infty}^2\right)d\tau+C\nonumber\\
&&\leq C(T),
\end{eqnarray}
which gives the desired estimate (\ref{2.66}). This completes the proof of Lemma 2.10.

From Lemma 2.2, (\ref{2.48}) and (\ref{2.66}), we have the following corollary.
\begin{Corollary}[Lower and upper  bounds of $v$ and $\theta$ for the case (iv) of Theorem 1.1]
Under the assumptions of Proposition 2.2, if the constants $\alpha,\beta$ satisfy the condition (iv) of Theorem 1.1, then it holds
\begin{equation}\label{2.71}
C^{-1}_6\leq v(t,x)\leq C_6,\quad
C_{11}^{-1}(T)\leq\theta(t,x)\leq C_{17}(T),\quad \forall\,(t,x)\in[0, T_1]\times\mathbb{R}.
\end{equation}
where the constants $C_6, C_{11}^{-1}(T),C_{17}(T)$ are given in Lemma 2.2, (\ref{2.48}) and (\ref{2.66}), respectively.
\end{Corollary}

Next, we use Lemmas 2.6-2.7 to show the following:
\begin{Lemma}[Lower and upper  bounds of $v$ and $\theta$ for the cases (i), (ii), (iii) and (v) of Theorem 1.1]
Under the assumptions of Proposition 2.2, if the constants $\alpha,\beta$ satisfy the condition (i) or (ii), or (iii), or (v) of Theorem 1.1, then it holds
\begin{equation}\label{2.72}
C^{-1}_{18}(T)\leq v(t,x)\leq C_{18}(T),\quad
C_{19}^{-1}(T)\leq\theta(t,x)\leq C_{19}(T),\quad \forall\, (t,x)\in[0, T_1]\times\mathbb{R},
\end{equation}
where $C_{18}(T),C_{19}(T)$ are some positive constants depending only on $T, \alpha,\beta,\lambda,\underline{V}, \bar{V}, \underline{\Theta}, \bar{\Theta}$, $\|v_0-1\|_2$, $\|u_0\|_2$ and $\|\theta_0-1\|$.
\end{Lemma}
\noindent{\bf Proof.}~For the case (i): $\alpha=0, \beta=-2$,  we have $v(t,x)\geq C_6^{-1}$ for all $(t,x)\in[0,T_1]\times\mathbb{R}$ by Remark 2.1. This together with Lemma 2.7 and the assumption $\alpha=0$ implies that
\begin{equation}\label{2.73}
\theta(t,x)\geq C^{-1}_{20}(T),\quad \forall\, (t,x)\in[0, T_1]\times\mathbb{R}
\end{equation}
for some positive constant $C_{20}(T)$ depending only on $T, \underline{V}, \bar{V}, \underline{\Theta}, \bar{\Theta}$ and  $\|(v_0-1, v_{0x}, u_0, \theta_0-1)\|$. Then it follows from (\ref{2.33}), (\ref{2.73}) and the assumption $\lambda\geq1$ that
$\|v\|_{L^\infty_{T_1,x}}\leq C_{21}(T)$
for some positive constant $C_{21}(T)$ depending only on $T, \alpha,\beta,\lambda,\underline{V}, \bar{V}, \underline{\Theta}, \bar{\Theta}$ and  $\|(v_0-1, v_{0x}, u_0, \theta_0-1)\|$. Consequently, we have
\begin{equation}\label{2.74}
C_6^{-1}\leq v(t,x)\leq C_{21}(T), \quad \forall\, (t,x)\in[0, T_1]\times\mathbb{R}.
 \end{equation}
With (\ref{2.73}) and (\ref{2.74}) in hand, some similar results as Corollary 2.2, Lemmas 2.8 and 2.9 hold. Then similar to Lemma 2.10, we can obtain
\begin{equation}\label{2.75}
\left\|\theta\right\|_{L^{\infty}_{T_1,x}}\leq C_{22}(T)
\end{equation}
for some  constant $C_{22}(T)>0$  depending only on  $T, \alpha,\beta,\lambda,\underline{V}, \bar{V}, \underline{\Theta}, \bar{\Theta}$, $\|v_0-1\|_2$, $\|u_0\|_2$ and $\|\theta_0-1\|$.

For the cases (ii): $\alpha<-2\beta-4, \beta\geq-\frac{3}{2}$  and (iii):  $\alpha<-\beta-\frac{5}{2}, -2\leq\beta<-\frac{3}{2}$, we have $\alpha<-\frac{1}{2}$ and $v(t,x)\geq C_6^{-1}$ for all $(t,x)\in[0,T_1]\times\mathbb{R}$ by Remark 2.1. Then we deduce
from Lemma 2.7 that
\begin{equation}\label{2.76}
\theta(t,x)\geq C^{-1}_{23}(T),\quad \forall\, (t,x)\in[0, T_1]\times\mathbb{R}
\end{equation}
for some positive constant $C_{23}(T)$ depending only on $T, \underline{V}, \bar{V}, \underline{\Theta}, \bar{\Theta}$ and  $\|(v_0-1, v_{0x}, u_0, \theta_0-1)\|$. (\ref{2.76}) together with (\ref{2.34}) and the assumption $\lambda\geq1$ yields  that $\|v\|_{L^\infty_{T_1,x}}\leq C_{24}(T)$
for some positive constant $C_{24}(T)$ depending only on $T, \alpha,\beta,\lambda,\underline{V}, \bar{V}, \underline{\Theta}, \bar{\Theta}$ and  $\|(v_0-1, v_{0x}, u_0, \theta_0-1)\|$. Then we can prove some similar results as  Corollary 2.2 and Lemmas 2.8 -2.10. Consequently, we obtain
\begin{equation}\label{2.77}
\left\|\theta\right\|_{L^{\infty}_{T_1,x}}\leq C_{25}(T)
\end{equation}
for some  constant $C_{25}(T)>0$  depending only on  $T, \alpha,\beta,\lambda,\underline{V}, \bar{V}, \underline{\Theta}, \bar{\Theta}$,  $\|v_0-1\|_2$, $\|u_0\|_2$ and $\|\theta_0-1\|$.

Finally, for the case (v): $\alpha>-2\beta-5, \beta<-3$, it follows  that $\alpha>1$ and  $v(t,x)\leq C_6$ for all $(t,x)\in[0,T_1]\times\mathbb{R}$ by Remark 2.1.  Then in view of  Lemma 2.7, we conclude
\begin{equation}\label{2.78}
\theta(t,x)\geq C^{-1}_{26}(T),\quad \forall\, (t,x)\in[0, T_1]\times\mathbb{R}
\end{equation}
for some positive constant $C_{26}(T)$ depending only on $T, \underline{V}, \bar{V}, \underline{\Theta}, \bar{\Theta}$ and  $\|(v_0-1, v_{0x}, u_0, \theta_0-1)\|$. It follows from (\ref{2.78}), (\ref{2.35}) and the assumption $\lambda\geq1$  that $\|\frac{1}{v}\|_{L^\infty_{T_1,x}}\leq C_{27}(T)$
for some positive constant $C_{27}(T)$ depending only on $T, \alpha,\beta,\lambda,\underline{V}, \bar{V}, \underline{\Theta}, \bar{\Theta}$ and  $\|(v_0-1, v_{0x}, u_0, \theta_0-1)\|$. Then we can show some similar results as  Corollary 2.2 and Lemmas 2.8-2.10. Thus we can obtain
\begin{equation}\label{2.79}
\left\|\theta\right\|_{L^{\infty}_{T_1,x}}\leq C_{28}(T)
\end{equation}
for some  constant $C_{28}(T)>0$  depending only on  $T, \alpha,\beta,\lambda,\underline{V}, \bar{V}, \underline{\Theta}, \bar{\Theta}$, $\|v_0-1\|_2$, $\|u_0\|_2$ and $\|\theta_0-1\|$.

 Letting $C_{18}(T)=\max\{C_6,C_{21}(T),C_{24}(T),C_{27}(T)\}$ and $C_{19}(T)=\max\{C_{20}(T), C_{22}(T), C_{23}(T),\\ C_{25}(T), C_{26}(T),C_{28}(T)\}$, then we complete the proof of Lemma 2.11.

From the proof of Lemma 2.11, we see that some similar results as Corollary 2.2 and Lemmas 2.8-2.9 still hold  for the cases (i), (ii), (iii) and (v) of Theorem 1.1.  Thus by combining Corollary 2.2 and Lemmas 2.8-2.9,  we have the following:
\begin{Corollary}
Under the assumptions of Proposition 2.2, there exists a  positive constant $C_{29}(T)$ depending only on $T, \alpha,\beta,\underline{V}, \bar{V}, \underline{\Theta}, \bar{\Theta}$, $\|v_0-1\|_2$, $\|u_0\|_2$ and $\|\theta_0-1\|_1$ such that for $0\leq t\leq T_1$,
\begin{equation}\label{2.80}
\|(v-1)(t)\|_2^2+\|u(t)\|_1^2+\|(\theta-1)(t)\|^2+\int_0^t\|(v_x,v_{xx},u_x,u_{xx},\theta_x)(\tau)\|^2d\tau\leq C_{29}(T).
\end{equation}
\end{Corollary}

For the estimate on $\|\theta_x(t)\|$, we have
\begin{Lemma}
Under the assumptions of Proposition 2.2, there exists a  positive constant $C_{30}(T)$ depending only on $T, \alpha,\beta,\underline{V}, \bar{V}, \underline{\Theta}, \bar{\Theta}$, $\|v_0-1\|_2$, $\|u_0\|_1$ and $\|\theta_0-1\|_1$ such that for $0\leq t\leq T_1$,
\begin{equation}\label{2.81}
\|\theta_x(t)\|+\int_0^t\|\theta_{xx}(\tau)\|^2d\tau\leq C_{30}(T).
\end{equation}
\end{Lemma}
\noindent{\bf Proof.}~Multiplying (\ref{2.78}) by $\theta^\lambda$,
we have
\begin{equation}\label{2.82}
C_vG(\theta)_t+p(v,\theta)u_x\theta^\lambda=\left(\frac{G(\theta)_x}{v}\right)_x\theta^\lambda+\frac{u^2_x}{v^{\alpha+1}}\theta^\lambda
\end{equation}
with $G(\theta)=\theta^{\lambda+1}$. Differentiating (\ref{2.82}) once with respect to $x$ and multiplying the resultant
equation by $G(\theta)_x$, and then integrating the equation over $[0,t]\times\mathbb{R}$ gives
\begin{equation}\label{2.83}
\aligned
&\frac{C_v}{2}\int_{\mathbb{R}}\left[G(\theta)_x\right]^2dx+\underbrace{\int_0^t\int_{\mathbb{R}}\left(p(v,\theta)u_x\theta^\lambda\right)_xG(\theta)_x dxd\tau}_{I_{11}}\\
=&\frac{C_v}{2}\int_{\mathbb{R}}\left[G(\theta_0)_x\right]^2dx
+\underbrace{\int_0^t\int_{\mathbb{R}}\left[\theta^\lambda\left(\frac{G(\theta)_x}{v}\right)_x\right]_xG(\theta)_xdxd\tau}_{I_{12}}+\underbrace{\int_0^t\int_{\mathbb{R}}\left[\theta^\lambda\frac{u^2_x}{v^{\alpha+1}}\right]_xG(\theta)_xdxd\tau}_{I_{13}}
\endaligned
\end{equation}
Using integration by parts,  the Cauchy inequality and (\ref{2.71})-(\ref{2.72}),  we have
\begin{equation}\label{2.84}
\aligned
&I_{11}=-\int_0^t\int_{\mathbb{R}}p(v,\theta)\theta^\lambda\left[\left(\frac{G(\theta)_x}{v}\right)_xv+\frac{G(\theta)_x}{v}v_x\right]dxd\tau
\\&\quad\leq\eta\int_0^t\int_{\mathbb{R}}v\theta^\lambda\left[\left(\frac{G(\theta)_x}{v}\right)_x\right]^2dxd\tau+C_\eta\int_0^t\int_{\mathbb{R}}
\left(|u_x|^2+|\theta_xv_x|^2\right)dxd\tau,
\endaligned
\end{equation}
\begin{equation}\label{2.85}
\aligned
&I_{12}=-\int_0^t\int_{\mathbb{R}}\theta^\lambda\left(\frac{G(\theta)_x}{v}\right)_x\left[\left(\frac{G(\theta)_x}{v}\right)_xv+\frac{G(\theta)_x}{v}v_x\right]dxd\tau
\\&\quad=-\int_0^t\int_{\mathbb{R}}v\theta^\lambda\left[\left(\frac{G(\theta)_x}{v}\right)_x\right]^2dxd\tau-\int_0^t\int_{\mathbb{R}}\theta^\lambda\frac{G(\theta)_xv_x}{v}
\left(\frac{G(\theta)_x}{v}\right)_xdxd\tau\\
&\quad\leq-\int_0^t\int_{\mathbb{R}}v\theta^\lambda\left[\left(\frac{G(\theta)_x}{v}\right)_x\right]^2dxd\tau+C_\eta\int_0^t\int_{\mathbb{R}}
|\theta_xv_x|^2dxd\tau,
\endaligned
\end{equation}
\begin{equation}\label{2.86}
\aligned
&I_{13}=-\int_0^t\int_{\mathbb{R}}\theta^\lambda\frac{u^2_x}{v^{\alpha+1}}\left[\left(\frac{G(\theta)_x}{v}\right)_xv+\frac{G(\theta)_x}{v}v_x\right]dxd\tau
\\
&\quad\leq\eta\int_0^t\int_{\mathbb{R}}v\theta^\lambda\left[\left(\frac{G(\theta)_x}{v}\right)_x\right]^2dxd\tau+C_\eta\int_0^t\int_{\mathbb{R}}
\left(u_x^4+|\theta_xv_x|^2\right)dxd\tau.
\endaligned
\end{equation}
Putting (\ref{2.83})-(\ref{2.85}) into (\ref{2.82}), we get by the smallness of $\eta$ that
\begin{equation}\label{2.87}
\aligned
&\int_{\mathbb{R}}\left[G(\theta)_x\right]^2dx+\int_0^t\int_{\mathbb{R}}v\theta^\lambda\left[\left(\frac{G(\theta)_x}{v}\right)_x\right]^2dxd\tau\\
=&C\int_{\mathbb{R}}\left[G(\theta_0)_x\right]^2dx
+C\int_0^t\int_{\mathbb{R}}\left(u_x^2+u_x^4+|\theta_xv_x|^2\right)dxd\tau.
\endaligned
\end{equation}
It follows from Corollary 2.4 and the Sobolev inequality that
\begin{equation}\label{2.88}
\aligned
&\int_0^t\int_{\mathbb{R}}\left(u_x^2+u_x^4+|\theta_xv_x|^2\right)dxd\tau\\
&\leq C(T)+\int_0^t\left(\sup_{0\leq\tau\leq t}\{\|u_x(\tau)\|^2\}\|u_x(\tau)\|\|u_{xx}(\tau)\|
+\sup_{0\leq\tau\leq t}\{\|v_x(\tau)\|_1^2\}\|\theta_x(\tau)\|^2\right)dt\\
&\leq C(T)+C(T)\int_0^t\left(\|u_x(\tau)\|_1^2
+\|\theta_x(\tau)\|^2\right)dt\leq C(T).
\endaligned
\end{equation}
Combining (\ref{2.87})-(\ref{2.88}), then it follows from Corollary 2.3 and Lemma 2.11 that
\begin{equation}\label{2.89}
\|\theta_x(t)\|^2+\int_0^t\int_{\mathbb{R}}\left[\left(\frac{G(\theta)_x}{v}\right)_x\right]^2dxd\tau\leq C(T)
\end{equation}
Moreover, we deduce  from (\ref{2.89}), (\ref{2.88}), (\ref{2.80}), the Cauchy inequality and the Sobolev inequality  that
\begin{eqnarray}\label{2.90}
&&\|\theta_x(t)\|^2+\int_0^t\|\theta_{xx}(\tau)\|^2d\tau\nonumber\\
&&\leq C(T)\left\{1+\int_0^t\int_{\mathbb{R}}\left(|\theta_x|^4+|\theta_xv_x|^2+|\theta_{xx}\theta_{x}^2|+|\theta_{xx}\theta_{x}v_x|\right)dt\right\}\nonumber\\
&&\leq \eta\int_0^t\|\theta_{xx}(\tau)\|^2d\tau+C(\eta,T)\int_0^t\int_{\mathbb{R}}\left(|\theta_x|^4+|\theta_xv_x|^2\right)dxdt+C(T)\nonumber\\
&&\leq \eta\int_0^t\|\theta_{xx}(\tau)\|^2d\tau
+C(\eta,T)\left(\int_0^t\|\theta_x(\tau)\|^3\|\theta_{xx}(\tau)\|dt+1\right)\nonumber\\
&&\leq 2\eta\int_0^t\|\theta_{xx}(\tau)\|^2d\tau
+C(\eta,T)\left(\int_0^t\sup_{0\leq\tau\leq t}\{\|\theta_x(\tau)\|^4\|\}\|\theta_x(\tau)\|^2dt+1\right)\nonumber\\
&&\leq 2\eta\int_0^t\|\theta_{xx}(\tau)\|^2d\tau
+C(\eta,T),
\end{eqnarray}
which, together with the smallness of $\eta$ leads to (\ref{2.81}) immediately. This completes the proof of Lemma 2.12.

Similarly, we can also obtain
\begin{Lemma}
Under the assumptions of Proposition 2.2, there exists a  positive constant $C_{31}(T)$ depending only on $T, \alpha,\beta,\underline{V}, \bar{V}, \underline{\Theta}, \bar{\Theta}$, $\|v_0-1\|_4$, $\|u_0\|_3$ and $\|\theta_0-1\|_3$ such that for $0\leq t\leq T_1$,
\begin{equation}\label{2.91}
\|v_{xx}(t)\|_2^2+\|(u_{xx},\theta_{xx})(t)\|_1^2+\int_0^t\left(\|v_{xx}(\tau)\|_3^2+\|(u_{xxx},\theta_{xxx})(\tau)\|_1^2\right)d\tau\leq C_{31}(T).
\end{equation}
\end{Lemma}

\noindent{\bf Proof of Proposition 2.2.} Proposition 2.2 follows immediately  from  Corollary 2.4 and Lemmas 2.12-2.13.

\section{Proof of Theorem 1.2}
\setcounter{equation}{0}
This section is devoted to proving Theorem 1.2. First, the local existence of the Cauchy problem (\ref{1.4})-(\ref{1.5}) has been given in Proposition 2.1. Suppose the local solution $(v,u, \theta)(t,x)$ has been extended to time $t=T_1\geq t_1$($t_1$ is given in Proposition 2.1), and satisfies the following a priori assumption:
\begin{equation}\label{3.1}
\frac{1}{2}\underline{\Theta}\leq\theta(t,x)\leq2\bar{\Theta},\quad \forall\,(t,x)\in[0,T_1]\times\mathbb{R}.
\end{equation}
Here and to the end of this subsection,  we assume that $0<\underline{\Theta}<1<\bar{\Theta}$ without  loss of generality.
Under the priori assumption (\ref{3.1}), if we can show that the specific volume $v(t,x)$ has uniform-in-time lower and upper bounds which are independent of $\gamma-1$, and some uniform-in-time, $\gamma-1$ independent energy estimates on $(v,u, \theta)(t,x)$, then we can close (\ref{3.1})  by using the smallness of $\gamma-1$. Then as usual, Theorem 1.2 can be obtained by combining the local existence result and the energy estimates. More precisely, we will prove the following:
\begin{Proposition} [A priori estimates]
 Under the assumptions of Theorem 1.2, suppose that $(v-1,u, \theta-1)(t,x)\in X(0, T_1; V_0, V_1; \frac{\underline{\Theta}}{2}, 2\bar{\Theta})$ is a solution of the Cauchy problem (\ref{1.4})-(\ref{1.5}) for some  positive constants $T_1, V_0$ and $V_1$.  Then there exist a  constant $C_4>0$ depending only on $\lambda, \underline{V}, \bar{V}, \underline{\Theta}, \bar{\Theta}$ and $\|(v_0-1, v_{0x}, u_0, \frac{\theta_0-1}{\sqrt{\gamma-1}})\|$ and a constant $C_5>0$ depending only on $\lambda, \underline{V}, \bar{V}, \bar{\Theta}, \underline{\Theta},\|v_0-1\|_4$ and $\left\|\left(u_0, \frac{\theta_0-1}{\sqrt{\gamma-1}}\right)\right\|_3$
 such that the following two estimates hold:
\begin{equation}\label{3.2}
C_4^{-1}\leq v(t,x)\leq C_4,\quad \underline{\Theta} \leq \theta(t,x)\leq \bar{\Theta},\quad\forall\,(t,x)\in[0,T_1]\times\mathbb{R},
\end{equation}
\begin{equation}\label{3.3}
\aligned
&\|(v-1)(t)\|^2_4+\left\|\left(u, \frac{\theta-1}{\sqrt{\gamma-1}}\right)(t)\right\|^2_3+\int_0^t\left(\|v_x(s)\|^2_4+\|(u_x,\theta_x)(s)\|^2_3\right)ds
\\
&\leq C_5\bigg(\|(v_0-1)(t)\|^2_4+\left\|\left(u_0, \frac{\theta_0-1}{\sqrt{\gamma-1}}\right)\right\|^2_3\bigg), \quad\forall\,t\in[0,T_1].
\endaligned
\end{equation}
\end{Proposition}

Proposition 3.1 can be obtained by a series of lemmas below. First, notice that Lemmas 2.1-2.2 and Remark 2.1 still hold under the assumptions of Proposition 3.1,  thus by combining Lemma 2.1 and the a priori assumption (\ref{3.1}),  we have the following lemma.
\begin{Lemma}[Basic energy estimates]
 Under the assumptions of Proposition 3.1, there  exists a  positive constant $C_{32}$ depending only on
 $\lambda,\underline{V}, \bar{V}, \underline{\Theta}$ and $\bar{\Theta}$  such that
\begin{equation}\label{3.4}
\aligned
 &\displaystyle\int_{\mathbb{R}}\left(R\phi(v)+\frac{u^2}{2}+\frac{R}{\gamma-1}\phi(\theta)\right)dx+\int_{\mathbb{R}}\frac{v_x^2}{v^{\beta+5}}dx
 +\displaystyle\int_0^t\int_{\mathbb{R}}\left(\frac{\theta_x^2}{v}+\frac{u_x^2}{v^{\alpha+1}}\right) dxd\tau\\
 &\leq C_{32}\left\|\left(v_0-1,u_0,\frac{\theta_0-1}{\sqrt{\gamma-1}},v_{0x}\right)\right\|^2
 \endaligned
 \end{equation}
for all $t\in[0, T_1]$, where the function $\phi(x):=x-1-\ln x$.
 \end{Lemma}

Similar to Lemma 2.2, we have
\begin{Lemma} Under the assumptions of Proposition 3.1, there  exists a  positive constant $C_{33}$ depending only on $ \lambda, \underline{V}, \bar{V}, \underline{\Theta}, \bar{\Theta}$ and $\|(v_0-1,v_{0x}, u_0, \frac{\theta_0-1}{\sqrt{\gamma-1}})\|$ such that
\begin{equation}\label{3.5}
C^{-1}_{33}\leq v(t,x)\leq C_{33}
 \end{equation}
for all $(t,x)\in[0, T_1]\times\mathbb{R}$.
\end{Lemma}
\noindent{\bf Proof.}~For the case (a) of Theorem 1.2, i.e., $\alpha=0,\beta=-2$,  we have $g(0,-2)=0$, where the function $g(\alpha,\beta)$ is defined in (\ref{2.15}). Then it follows from (\ref{2.16}) for $\alpha=0$ and $\beta=-2$ that
\begin{eqnarray}\label{3.6}
&&\displaystyle\int_\mathbb{R}\frac{ v_x^2}{v^2}dx+\displaystyle\int_0^t\int_\mathbb{R}\frac{ R\theta v_x^2}{v}dxd\tau+\int_0^t\int_\mathbb{R}\left[\left(\frac{v_x}{v^{2}}\right)_x\right]^2dxd\tau
\nonumber\\
&&\leq C\displaystyle\left(\int_\mathbb{R}\frac{v_{0x}^2}{v_0^{2}}dx+\|u(t)\|^2+\|u_0\|^2\right)
+\int_0^t\int_\mathbb{R}\frac{u_x^2}{v}dxd\tau\displaystyle+\int_0^t\int_\mathbb{R}\frac{R\theta_xv_x}{v^{\alpha+2}}dxd\tau
,
\end{eqnarray}
Utilizing the Cauchy inequality, the a priori assumption (\ref{3.1}) and Lemma 3.1 with $\alpha=0$,  we infer
\begin{equation}\label{3.7}
\aligned
\int_0^t\int_\mathbb{R}\frac{R\theta_xv_x}{v^{2}}dxd\tau&\leq\displaystyle\frac{1}{2}\int_0^t\int_\mathbb{R}\frac{R\theta v_x^2}{v}dxd\tau+C\int_0^t\int_\mathbb{R}\frac{\theta_x^2}{\theta v}dxd\tau\\
&\leq\displaystyle\frac{1}{2}\int_0^t\int_\mathbb{R}\frac{R\theta v_x^2}{v}dxd\tau+C(\underline{\Theta})\int_0^t\int_\mathbb{R}\frac{\theta_x^2}{v}dxd\tau\\
&\leq\displaystyle\frac{1}{2}\int_0^t\int_\mathbb{R}\frac{R\theta v_x^2}{v}dxd\tau+C(\lambda,\underline{V}, \bar{V}, \underline{\Theta}, \bar{\Theta})\left\|\left(v_0-1,u_0,\frac{\theta_0-1}{\sqrt{\gamma-1}},v_{0x}\right)\right\|^2.
\endaligned
\end{equation}
Inserting  (\ref{3.7}) into (\ref{3.6}), we get from (\ref{3.1}) and Lemma 3.1 that
\begin{eqnarray}\label{3.8}
&&\displaystyle\int_\mathbb{R}\frac{ v_x^2}{v^2}dx+\displaystyle\int_0^t\int_\mathbb{R}\left(\frac{ v_x^2}{v}+\left[\left(\frac{v_x}{v^{2}}\right)_x\right]^2\right)dxd\tau
\leq C(\lambda,\underline{V}, \bar{V}, \underline{\Theta}, \bar{\Theta})\left\|\left(v_0-1,u_0,\frac{\theta_0-1}{\sqrt{\gamma-1}},v_{0x}\right)\right\|^2.
\end{eqnarray}
By Remark 2.1, $v(t,x)\geq C_6^{-1}$ for the case (a) of Theorem 2.1. Now we prove the upper bound of $v(t,x)$. Similar to (\ref{2.38}), we derive  from Lemma 3.1 and (\ref{3.8}) that
\begin{equation}\label{3.9}
\aligned
v^{\frac{1}{2}}
&\leq C+C\left\|\sqrt{\phi(v)}\right\|\left\|\frac{v_x}{v}\right\|\\
&\leq C+C(\lambda,\underline{V}, \bar{V}, \underline{\Theta}, \bar{\Theta})\left\|\left(v_0-1,u_0,\frac{\theta_0-1}{\sqrt{\gamma-1}},v_{0x}\right)\right\|^2,
\endaligned
\end{equation}
which implies that
\begin{equation}\label{3.10}
\left\|v\right\|_{L^\infty_{T_1,x}}
\leq C(\lambda,\underline{V}, \bar{V}, \underline{\Theta}, \bar{\Theta})\left\|\left(v_0-1,u_0,\frac{\theta_0-1}{\sqrt{\gamma-1}},v_{0x}\right)\right\|^4:=C_{34}.
\end{equation}

For the case  (b) of Theorem 1.2, i.e., $\beta=2\alpha-3,-3\leq\beta<-2$,  we have Lemma 2.2 holds.  Therefore, we finish the proof of Lemma 3.2 by letting $C_{33}=\max\{C_6, C_{34}\}$.

As a consequence of Lemmas 3.1 and 3.2, we have
\begin{Corollary} Under the assumptions of Proposition 3.1, there  exists a  positive constant $C_{35}$ depending only on $ \lambda, \underline{V}, \bar{V}, \underline{\Theta}, \bar{\Theta}$ and $\|(v_0-1,v_{0x}, u_0, \frac{\theta_0-1}{\sqrt{\gamma-1}})\|$ such that
\begin{equation}\label{3.11}
 \displaystyle\left\|\left(v-1,u,\frac{\theta-1}{\sqrt{\gamma-1}},v_x\right)(t)\right\|^2+\int_0^t\|(u_x,\theta_x)(\tau)\|^2 d\tau\leq C_{35}\left\|\left(v_0-1,u_0,\frac{\theta_0-1}{\sqrt{\gamma-1}},v_{0x}\right)\right\|^2
 \end{equation}
for all $(t,x)\in[0, T_1]\times\mathbb{R}$.
\end{Corollary}

Next, we deduce an estimate on $\displaystyle\int_0^t\|v_x(\tau)\|_1^2d\tau$.
\begin{Lemma} Under the assumptions of Proposition 3.1, there  exists a  positive constant $C_{36}$ depending only on $ \lambda, \underline{V}, \bar{V}, \underline{\Theta}, \bar{\Theta}$ and $\|(v_0-1,v_{0x}, u_0, \frac{\theta_0-1}{\sqrt{\gamma-1}})\|$   such that
\begin{equation}\label{3.12}
 \displaystyle\left\|v_x(t)\right\|^2+\int_0^t\|v_x(\tau)\|_1^2 d\tau\leq C_{36}\left\|\left(v_0-1,u_0,\frac{\theta_0-1}{\sqrt{\gamma-1}},v_{0x}\right)\right\|^6
 \end{equation}
for all $(t,x)\in[0, T_1]\times\mathbb{R}$.
\end{Lemma}
\noindent{\bf Proof.}~If the constants $\alpha,\beta$ satisfy the condition (a): $\alpha=0,\beta=-2$, then we have from (\ref{3.8}) and Lemma 3.2 that
\begin{equation}\label{3.13}
 \displaystyle\left\|v_x(t)\right\|^2+\int_0^t\|v_x(\tau)\|^2 d\tau\leq C(\lambda,\underline{V}, \bar{V}, \underline{\Theta}, \bar{\Theta})\left\|\left(v_0-1,u_0,\frac{\theta_0-1}{\sqrt{\gamma-1}},v_{0x}\right)\right\|^2.
 \end{equation}
On the other hand,  (\ref{3.8}) and Lemma 3.2 also imply that
\begin{equation}\label{3.14}
\aligned
 &\displaystyle\left\|v_x(t)\right\|^2+\int_0^t\|v_x(\tau)\|_1^2 d\tau\\
 &\leq C(\lambda,\underline{V}, \bar{V}, \underline{\Theta}, \bar{\Theta})\left(\left\|\left(v_0-1,u_0,\frac{\theta_0-1}{\sqrt{\gamma-1}},v_{0x}\right)\right\|^2+\int_0^t\int_\mathbb{R}\left(v_x^4+|v_{xx}v_x^2|\right)dxd\tau\right).
 \endaligned
 \end{equation}
 It follows from the Cauchy inequality and (\ref{3.13}) that
\begin{equation}\label{3.15}
\aligned
 \displaystyle\int_0^t\int_\mathbb{R}\left(v_x^4+|v_{xx}v_x^2|\right)dxd\tau&\leq\eta\int_0^t\|v_{xx}(\tau)\|^2d\tau+C_\eta\displaystyle\int_0^t\int_\mathbb{R}v_x^4dxd\tau\\
 &\leq\frac{\eta}{2}\int_0^t\|v_{xx}(\tau)\|^2d\tau+C_\eta\displaystyle\int_0^t\int_\mathbb{R}\|v_x(\tau)\|^3\|v_{xx}(\tau)\|dxd\tau\\
  &\leq\eta\int_0^t\|v_{xx}(\tau)\|^2d\tau+C_\eta\displaystyle\int_0^t\sup_{0\leq\tau\leq t}\{\|v_x(\tau)\|^4\}\|v_{x}(\tau)\|^2dxd\tau\\
  &\leq\eta\int_0^t\|v_{xx}(\tau)\|^2d\tau+C_\eta\left\|\left(v_0-1,u_0,\frac{\theta_0-1}{\sqrt{\gamma-1}},v_{0x}\right)\right\|^6.
 \endaligned
 \end{equation}
Combining (\ref{3.14})-(\ref{3.15}), and using the smallness of $\eta$, we have
\begin{equation}\label{3.16}
 \displaystyle\left\|v_x(t)\right\|^2+\int_0^t\|v_x(\tau)\|_1^2 d\tau\leq C_{37}\left\|\left(v_0-1,u_0,\frac{\theta_0-1}{\sqrt{\gamma-1}},v_{0x}\right)\right\|^6,
 \end{equation}
where $C_{37}$ is a positive constant depending only on $ \lambda, \underline{V}, \bar{V}, \underline{\Theta}, \bar{\Theta}$ and $\|(v_0-1,v_{0x}, u_0, \theta_0-1)\|$ .

If $\alpha,\beta$ satisfy the condition (b): $\beta=2\alpha-3$, $-3\leq\beta<-2$, then we get by integrating (\ref{2.13}) with respect to $t$ and $x$ over $[0,t]\times\mathbb{R}$ that
\begin{eqnarray}\label{3.17}
&&\displaystyle\frac{1}{4}\int_\mathbb{R}\frac{v_x^2}{v^{2\alpha+2}}dx+\displaystyle\int_0^t\int_\mathbb{R}\frac{R\theta v_x^2}{v^{\alpha+1}}dxd\tau+\int_0^t\int_\mathbb{R}\left[\left(\frac{v_x}{v^{\frac{\alpha+\beta+6}{2}}}\right)_x\right]^2dxd\tau
\nonumber\\
&&\leq C\displaystyle\left(\int_\mathbb{R}\frac{v_{0x}^2}{v_0^{2\alpha+2}}dx+\|u(t)\|^2+\|u_0\|^2\right)
+\int_0^t\int_\mathbb{R}\frac{u_x^2}{v^{\alpha+1}}dxd\tau\nonumber\\
&&\quad\displaystyle+\int_0^t\int_\mathbb{R}\left\{\frac{v_{xx}}{v^{\beta+5}}-\frac{\beta+5}{2}\frac{v_x^2}{v^{\beta+6}}\right\}_x
\frac{v_x^2}{v^{\alpha+1}}dxd\tau,
\end{eqnarray}
where we have used the fact that
\[
\int_{\mathbb{R}}\frac{u v_x}{v^{\alpha+1}}\,dx\leq\frac{1}{4}\int_{\mathbb{R}}\frac{v_x^2}{v^{2\alpha+2}}\,dx+C\int_{\mathbb{R}}u^2\,dx.
\]
Using integration by parts, we obtain
\begin{equation}\label{3.18}
\aligned
 &\displaystyle\int_0^t\int_\mathbb{R}\left\{\frac{v_{xx}}{v^{\beta+5}}-\frac{\beta+5}{2}\frac{v_x^2}{v^{\beta+6}}\right\}_x
\frac{v_x}{v^{\alpha+1}}dxd\tau\\
 =&-\displaystyle\int_0^t\int_\mathbb{R}\left\{\frac{1}{v^{\beta+2}}\left(\frac{-v_{x}}{v^{3}}\right)_x+\frac{\beta-1}{2}\frac{v_x^2}{v^{\beta+6}}\right\}_x
 \frac{v_x}{v^{\alpha+1}}dxd\tau\\
 =&-\displaystyle\int_0^t\int_\mathbb{R}\frac{1}{v^{\beta+2}}\left(\frac{v_x}{v^{\alpha+1}}\cdot v^{\alpha-2}\right)_x
 \left(\frac{v_x}{v^{\alpha+1}}\right)_xdxd\tau+
\frac{\beta-1}{2}\int_0^t\int_\mathbb{R}\frac{v_x^2}{v^{\beta+6}}
 \left(\frac{v_x}{v^{\alpha+1}}\right)_xdxd\tau\\
 =&-\displaystyle\int_0^t\int_\mathbb{R}\frac{1}{v^{\beta-\alpha+4}}\left[\left(\frac{v_x}{v^{\alpha+1}}\right)_x\right]^2
dxd\tau+\left[\frac{\beta-1}{2}-(\alpha-2)\right]\displaystyle\int_0^t\int_\mathbb{R}\frac{v_x^2}{v^{\beta+6}}
 \left(\frac{v_x}{v^{\alpha+1}}\right)_xdxd\tau\\
 =&-\displaystyle\int_0^t\int_\mathbb{R}\frac{1}{v^{\beta-\alpha+4}}\left[\left(\frac{v_x}{v^{\alpha+1}}\right)_x\right]^2
dxd\tau.
 \endaligned
 \end{equation}
where in the last step of (\ref{3.18}), we have used the assumption that $\beta=2\alpha-3$.

Substituting (\ref{3.18}) into (\ref{3.17}), then it follows from Lemmas 3.1-3.2 and the a priori assumption (\ref{3.1}) that
\begin{equation}\label{3.19}
 \displaystyle\left\|v_x(t)\right\|^2+\int_0^t\|v_x(\tau)\|^2 d\tau+\int_0^t\int_\mathbb{R}\left[\left(\frac{v_x}{v^{\alpha+1}}\right)_x\right]^2dxd\tau
 \leq C\left\|\left(v_0-1,u_0,\frac{\theta_0-1}{\sqrt{\gamma-1}},v_{0x}\right)\right\|^2.
 \end{equation}
Moreover, we  deduce from (\ref{3.19}), (\ref{3.15}) and Lemma 3.2 that
\begin{equation}\label{3.20}
 \aligned
 \displaystyle\left\|v_x(t)\right\|^2+\int_0^t\|v_x(\tau)\|_1^2 d\tau
 &\leq C\left\|\left(v_0-1,u_0,\frac{\theta_0-1}{\sqrt{\gamma-1}},v_{0x}\right)\right\|^2+C\int_0^t\int_\mathbb{R}\left(v_x^4+|v_{xx}v_x^2|\right)dxd\tau\\
 &\leq\eta\int_0^t\|v_{xx}(\tau)\|^2d\tau+C_\eta\left\|\left(v_0-1,u_0,\frac{\theta_0-1}{\sqrt{\gamma-1}},v_{0x}\right)\right\|^6,
  \endaligned
 \end{equation}
which together with the smallness of $\eta$ implies that
 \begin{equation}\label{3.21}
 \aligned
 \displaystyle\left\|v_x(t)\right\|^2+\int_0^t\|v_x(\tau)\|_1^2 d\tau
 &\leq C_{38}\left\|\left(v_0-1,u_0,\frac{\theta_0-1}{\sqrt{\gamma-1}},v_{0x}\right)\right\|^6.
  \endaligned
 \end{equation}
 Here $C_{38}$ is a positive constant depending on $\lambda, \underline{V}, \bar{V}, \underline{\Theta}, \bar{\Theta}$ and $\|(v_0-1,v_{0x}, u_0, \frac{\theta_0-1}{\sqrt{\gamma-1}})\|$. Letting $C_{36}=\max\{C_{37},C_{38}\}$, then we can get  (\ref{3.12}) from (\ref{3.16}) and (\ref{3.21}). This completes the proof of Lemma 3.3.

For the estimate on $\|u_x(t)\|$, we have
\begin{Lemma}
Under the assumptions of Proposition 3.1,  there exits a constant $C_{39}>0$  depending only on  $\lambda,\underline{V}, \bar{V}, \underline{\Theta}, \bar{\Theta}$, $\|v_0-1\|_2$, $\|u_0\|_1$ and $\|\frac{\theta_0-1}{\sqrt{\gamma-1}}\|$  such that  for $0\leq t\leq T_1$,
\begin{equation}\label{3.22}
\left\|(u_x,v_{xx})(t)\right\|^2+\int_0^t\|u_{xx}(\tau)\|^2d\tau\leq C_{39}\left(\|v_0-1\|_2^2+\|u_0\|_1^2+\left\|\frac{\theta_0-1}{\sqrt{\gamma-1}}\right\|^2\right).
\end{equation}
\end{Lemma}
\noindent{\bf Proof.}~Integrating (\ref{2.56}) with respect to $t$ and $x$ over $[0,t]\times\mathbb{R}$, we have from Lemma  3.2 that
\begin{equation}\label{3.23}
 \displaystyle\|(u_x, v_{xx})(t)\|^2+\displaystyle\int_0^t\|u_{xx}(\tau)\|^2\,d\tau
 \leq C\left(\|(v_{0xx},u_{0x})\|^2+I_{14}+I_{15}\right),
 \end{equation}
where
\[
\begin{array}{rl}
 &I_{14}=\displaystyle\int_0^t\int_\mathbb{R}\left\{\left|v_{xx}^2u_x\right|+|\theta_xu_{xx}|\right\}dxd\tau,\\[3mm]
 &I_{15}=\displaystyle\int_0^t\int_\mathbb{R}\left\{\left|\theta v_xu_{xx}\right|+\left|v_xv_{xx}u_{xx}\right|+\left|u_xv_xu_{xx}\right|+\left|v_{x}^3u_{xx}\right|\right\}dxd\tau.
 \end{array}
\]
It follows from the Sobolev inequality,  the Cauchy inequality,  the a priori assumption (\ref{3.1}), and  Lemmas 3.1 and 3.3 that
\begin{equation}\label{3.24}
\aligned
 I_{14}&\leq \int_0^t\left\{\|u_{x}(\tau)\|^{\frac{1}{2}}\|u_{xx}(\tau)\|^{\frac{1}{2}}\|v_{xx}(\tau)\|^2
 +\|u_{xx}(\tau)\|\|\theta_{x}(\tau)\|\right\}\,d\tau\\
 &\leq\eta\int_0^t\|u_{xx}(\tau)\|^{2}d\tau+
 C_\eta\left\{\int_0^t\left(\|v_{xx}(\tau)\|^{4}+\|(u_x,\theta_x)(\tau)\|^2\right)d\tau\right\}
 \\
 &\leq\eta\int_0^t\|u_{xx}(\tau)\|^{2}d\tau+
 C_\eta\left\{\int_0^t\|v_{xx}(\tau)\|^{4}d\tau+\left\|\left(v_0-1,u_0,\frac{\theta_0-1}{\sqrt{\gamma-1}},v_{0x}\right)\right\|^2\right\},
\endaligned
\end{equation}
\begin{eqnarray}\label{3.25}
I_{15} &&\leq\eta\int_0^t\|u_{xx}(\tau)\|^{2}d\tau+C_\eta\int_0^t\int_{\mathbb{R}}\left(|u_xv_x|^2+|v_xv_{xx}|^2+|v_x|^6+|v_x|^2\right)dxd\tau\nonumber\\
 &&\leq\eta\int_0^t\|u_{xx}(\tau)\|^{2}d\tau+C_\eta\int_0^t\left(\|u_x(\tau)\|\|u_{xx}(\tau)\|\sup_{0\leq\tau\leq t}\{\|v_x(\tau)\|^2\}\right.\nonumber\\
 &&\left.\qquad+\|v_{x}(\tau)\|\|v_{xx}(\tau)\|^3
 +\sup_{0\leq\tau\leq t}\{\|v_{x}(\tau)\|^4\}\|v_{xx}(\tau)\|^2+\|v_{x}(\tau)\|^2\right)d\tau\\
 &&\leq2\eta\int_0^t\|u_{xx}(\tau)\|^{2}d\tau+C_\eta\left(\int_0^t\|v_{x}(\tau)\|_1^2\|v_{xx}(\tau)\|^2d\tau+\left\|\left(v_0-1,u_0,\frac{\theta_0-1}{\sqrt{\gamma-1}},v_{0x}\right)\right\|^{10}\right).\nonumber
\end{eqnarray}
Combining (\ref{3.23})-(\ref{3.25}), then (\ref{3.22}) follows from Gronwall's inequality and the smallness of $\eta$. This completes the proof of Lemma 3.4.

The following lemma gives an estimate on $\|\theta_x(t)\|$.
\begin{Lemma}
Under the assumptions of Proposition 3.1,  there exits a constant $C_{40}>0$  depending only on  $\lambda,\underline{V}, \bar{V}, \underline{\Theta}, \bar{\Theta}$, $\|v_0-1\|_2$, $\|u_0\|_1$ and $\|\frac{\theta_0-1}{\sqrt{\gamma-1}}\|_1$  such that  for $0\leq t\leq T_1$,
\begin{equation}\label{3.26}
\left\|\frac{\theta_x}{\sqrt{\gamma-1}}(t)\right\|^2+\int_0^t\|\theta_{xx}(\tau)\|^2d\tau\leq C_{40}\left(\|v_0-1\|_2^2+\|u_0\|_1^2+\left\|\frac{\theta_0-1}{\sqrt{\gamma-1}}\right\|_1^2\right).
\end{equation}
\end{Lemma}

Lemma 3.5 can be deduced by using the argument in the proof of Lemmas 2.12 and 3.4.  Since  the proof is similar to that of Lemma 2.12,  we omit it here for brevity.

Under the a priori assumption (\ref{3.1}), we have  obtained Corollary 3.1 and Lemmas 3.1-3.5.  Now we  begin to close the a priori assumption (\ref{3.1}). First, it follows from the Sobolev inequality, Corollary 3.1 and Lemma 3.5 that
\begin{equation}\label{3.27}
\|\theta(t)-1\|_{L^\infty_{x}}\leq\|\theta(t)-1\|^{\frac{1}{2}}\|\theta_x(t)\|^{\frac{1}{2}}
\leq C_{41}(\gamma-1)^{\frac{1}{2}}\left(\|v_0-1\|_2+\|u_0\|_1+\left\|\frac{\theta_0-1}{\sqrt{\gamma-1}}\right\|_1\right),
\end{equation}
where $C_{41}$ is a positive constant depending only on $\lambda,\underline{V}, \bar{V}, \underline{\Theta}, \bar{\Theta}$, $\|v_0-1\|_2$, $\|u_0\|_1$ and $\|\frac{\theta_0-1}{\sqrt{\gamma-1}}\|_1$.

Since $\theta=\frac{A}{R}v^{1-\gamma}\exp\left(\frac{\gamma-1}{R}s\right)$ and $\bar{s}:=\frac{R}{\gamma-1}\ln\frac{A}{R}$, we have
\[\aligned
\theta-1&=\frac{A}{R}v^{1-\gamma}\exp\left(\frac{\gamma-1}{R}s\right)-\frac{A}{R}\exp\left(\frac{\gamma-1}{R}\bar{s}\right)\\
&=\frac{A}{R}v^{1-\gamma}\exp\left(\frac{\gamma-1}{R}s\right)+
\frac{A}{R}\left(\exp\left(\frac{\gamma-1}{R}s\right)-\exp\left(\frac{\gamma-1}{R}\bar{s}\right)\right),
\endaligned\]

\[
\theta_{x}=\frac{A}{R}(1-\gamma)v^{-\gamma}v_x\exp\left(\frac{\gamma-1}{R}s\right)+
\frac{A}{R}v^{1-\gamma}\exp\left(\frac{\gamma-1}{R}s\right)\frac{\gamma-1}{R}s_x.
\]
Consequently,
\begin{equation}\label{3.28}
\aligned
&\left\|\frac{\theta_0-1}{\sqrt{\gamma-1}}\right\|\leq C(\gamma-1)^{\frac{1}{2}}\exp\left(\frac{\gamma-1}{R}\|s_0\|_{L^\infty}\right)\left(\left(\inf v_0(x)\right)^{-\gamma}\|v_0-1\|+\|s-\bar{s}\|\right),\\
&\left\|\frac{\theta_{0x}}{\sqrt{\gamma-1}}\right\|\leq C(\gamma-1)^{\frac{1}{2}}\exp\left(\frac{\gamma-1}{R}\|s_0\|_{L^\infty}\right)\left(\left(\inf v_0(x)\right)^{-\gamma}\|v_{0x}\|+\left(\inf v_0(x)\right)^{1-\gamma}\|s_{0x}\|\right).
\endaligned
\end{equation}
Since $\frac{\gamma-1}{R}\|s_0\|_{L^\infty}$ and $\inf v_0(x)$ are assumed to be bounded by some constants independent of $\gamma-1$, we deduce from (\ref{3.27})-(\ref{3.28}) that
\begin{equation}\label{3.29}
\|\theta-1\|_{L^\infty_{T_1,x}}\leq C_{42}(\gamma-1)^{\frac{1}{2}}\left(\|v_0-1\|_2+\|u_0\|_1\right)+C_{42}(\gamma-1)\left\|(v_0-1,s_0-\bar{s})\right\|_1,
\end{equation}
where $C_{42}$ is a positive constant depending only on $\lambda,\underline{V}, \bar{V}, \underline{\Theta}, \bar{\Theta}$,  $\|v_0-1\|_2$, $\|u_0\|_1$ and $\|\frac{\theta_0-1}{\sqrt{\gamma-1}}\|_1$.

Thus if $\gamma-1$ is sufficiently small such that
\begin{equation}\label{3.30} C_{42}(\gamma-1)^{\frac{1}{2}}\left(\|v_0-1\|_2+\|u_0\|_1\right)+C_{42}(\gamma-1)\left\|(v_0-1,s_0-\bar{s})\right\|_1<\min\{\bar{\Theta}-1,1-\underline{\Theta}\},
\end{equation}
then we have
\begin{equation}\label{3.31}
\theta(t,x)\leq\|\theta(t,x)-1\|_{L^\infty_{T_1,x}}+1\leq1+\min\{\bar{\Theta}-1,1-\underline{\Theta}\}\leq\bar{\Theta},
\end{equation}
and \begin{equation}\label{3.32}
\theta(t,x)\geq1-\|\theta(t,x)-1\|_{L^\infty_{T_1,x}}\geq1-\min\{\bar{\Theta}-1,1-\underline{\Theta}\}\leq\bar{\Theta}\geq\underline{\Theta}
\end{equation}
for all $(t,x)\in[0,T_1]\times\mathbb{R}$. This closes the priori assumption (\ref{3.1}).

Combining Corollary 3.1 and Lemmas 3.2-3.5, we arrive at
\begin{equation}\label{3.33}
\aligned
 &\displaystyle\|(v-1)(t)\|_2^2+\left\|\left(u,\frac{\theta-1}{\sqrt{\gamma-1}}\right)(t)\right\|_1^2+\int_0^t\left(\|v_x(\tau)\|_1^2+\|(u_x,\theta_x)(\tau)\|_1^2\right) d\tau\\
 &\leq C_{43}\left(\|v_0-1\|_2^2+\left\|\left(u_0,\frac{\theta_0-1}{\sqrt{\gamma-1}}\right)\right\|_1^2\right), \quad\forall t\in[0,T_1],
 \endaligned
 \end{equation}
where $C_{43}$ is a positive  constant depending only on $\lambda,\underline{V}, \bar{V}, \underline{\Theta}, \bar{\Theta}$, $\|v_0-1\|_2$, $\|u_0\|_1$ and $\|\frac{\theta_0-1}{\sqrt{\gamma-1}}\|_1$.

By adopting the same argument as above, we can also obtain
 \begin{equation}\label{3.34}
\aligned
 &\displaystyle\|v_{xx}(t)\|_2^2+\left\|\left(u_{xx},\frac{\theta_{xx}}{\sqrt{\gamma-1}}\right)(t)\right\|_1^2+\int_0^t\left(\|v_{xx}(\tau)\|_3^2+\|(u_{xxx},\theta_{xxx})(\tau)\|_1^2\right) d\tau\\
 &\leq C_{44}\left(\|v_0-1\|_4^2+\left\|\left(u_0,\frac{\theta_0-1}{\sqrt{\gamma-1}}\right)\right\|_3^2\right), \quad\forall t\in[0,T_1],
 \endaligned
 \end{equation}
where $C_{44}$ is a positive   constant depending only on $\lambda,\underline{V}, \bar{V}, \underline{\Theta}, \bar{\Theta}$,  $\|v_0-1\|_4$ and $\|(u_0, \frac{\theta_0-1}{\sqrt{\gamma-1}})\|_3$.

 Proposition 3.1 follows immediately from Lemma 3.2,  and (\ref{3.31})-(\ref{3.34}). And then Theorem 1.2 can be obtained  by the standard continuation argument. Hence, we finish the proof of Theorem 1.2.

%%%%%%%%%%%%%%%%%%%%%%%%%%%%%%%%%%%%%%%%%%%%%%%%%%%%%%%%%%%%%%%%%%%%%%%%%%%%%%%%%%%%%%%%%%%%%%%%%%%%%%%

\bigbreak

\begin{center}
{\bf Acknowledgement}
\end{center}
 ZZC was supported by the Tian Yuan Foundation of China under contract 11426031, the National Natural Science Foundation of China under contract 11501003, and the  Doctoral Scientific Research Funds of Anhui University under  contract J10113190005. HJZ was support by two grants from the National Natural Science Foundation of China under contracts 10925103 and 11271160 respectively. This work was also supported by a grant from the National Natural Science Foundation of China under contract 11261160485 and the Fundamental Research Funds for the Central Universities.

\end{document}